\documentclass[a4paper,11pt]{article} 
\usepackage[english]{babel} 
\usepackage[matrix,arrow,tips,curve]{xy}
\usepackage{amsmath,amssymb,amsthm}
\usepackage{enumerate}

\newcounter{subsubsec}

\input{psfig.tex}
\input{psfig.sty}

\newtheorem{teo}{Theorem}[subsubsec]
\newcommand{\pr}[1]{\ensuremath{\mathbb{P}^{#1}}}

\newcommand{\ldo}{.\,.\,.\,}
\newcommand{\Pic}{\operatorname{Pic}}

\newcommand{\ph}{\ensuremath{\varphi}}
\newcommand{\lis}[1]{\ensuremath{\{x_1,\ldo,x_{#1}}\}}
\newcommand{\som}[1]{\ensuremath{x_1+\cdots+x_{#1}}}
\newcommand{\fx}{\ensuremath{\Sigma_X}}
\newcommand{\fy}{\ensuremath{\Sigma_Y}}
\newcommand{\Z}{\ensuremath{\mathbb{Z}}}
\newcommand{\f}{\ensuremath{\Sigma}}
\newcommand{\Q}{\ensuremath{\mathbb{Q}}}
\newcommand{\RelInt}{\operatorname{Rel\,Int}}

\newcommand{\NE}{\operatorname{NE}}
\newcommand{\N}{\ensuremath{\mathcal{N}_1}}
\newcommand{\A}{\ensuremath{\mathcal{A}_1}}
\newcommand{\ld}{.\,.\,}

\newtheorem{cor}[teo]{Corollary}
\newtheorem{emdefi}[teo]{Definition}
\newtheorem{prop}[teo]{Proposition}
\newtheorem{lemma}[teo]{Lemma}

\newenvironment{dimo}[1][Proof.] {\begin{proof}[#1]} {\end{proof}}
\newcommand{\PC}{\operatorname{PC}}
\newenvironment{defi}{\begin{emdefi} \em}{\end{emdefi}}

\begin{document}

{\noindent\bfseries\Large Contractible classes in toric varieties }

\bigskip

{\noindent \bfseries\large Cinzia Casagrande }



\medskip 

{\small
\noindent Dipartimento di Matematica\\ Universit\`a di Roma ``La Sapienza''
\\ piazzale Aldo Moro, 5\\ 
00185 Roma \\
ccasagra@mat.uniroma1.it}

\bigskip

\bigskip

\noindent Let $X$ be a smooth, complete toric variety. 
Let $\A(X)$ be the group of algebraic 1-cycles on $X$ modulo
numerical equivalence and $\N(X)=\A(X)\otimes_{\Z}\Q\,$. Consider 
in $\N(X)$ the cone $\NE(X)$ generated by classes of curves on $X$. 
It is a well-known result due to M.~Reid~\cite{reid} that
$\NE(X)$ is closed, polyhedral and generated by classes of invariant
curves on $X$. The variety $X$ is projective if and only if $\NE(X)$
is strictly convex; in this case, a 1-dimensional face of $\NE(X)$ is
called an extremal ray. It is shown in~\cite{reid} that every extremal
ray admits a contraction to a projective toric variety.

We think of $\A(X)$ as a lattice in the $\Q\,$-vector space
$\N(X)$. Suppose that $X$ is projective. For every extremal ray
$R\subset \NE(X)$, we choose the primitive class in $R\cap\A(X)$; we
call this class an extremal class. The set $\mathcal{E}$
of extremal classes is a generating
set for the cone $\NE(X)$, namely $\NE(X)=\sum_{\gamma\in\mathcal{E}}
\Q_{\,\geq 0}\,\gamma$. For many purposes it would be useful to have a
linear decomposition with 
\emph{integral} coefficients: for instance, what can we say about curves
having minimal degree with respect to some ample line bundle on $X$?
It is an open question whether extremal classes generate
$\NE(X)\cap\A(X)$ as a semigroup. In this paper we introduce a set
$\mathcal{C}\supseteq\mathcal{E}$ of classes in $\NE(X)\cap\A(X)$
which is a set of generators of $\NE(X)\cap\A(X)$ as a semigroup.
Classes in $\mathcal{C}$ are geometrically characterized by
``contractibility'': 

\smallskip

\noindent \textbf{Definition \ref{pioggia}.}
Let $\gamma\in\NE(X)\cap\A(X)$ be primitive along $\A(X)\cap\Q_{\,\geq
  0}\,\gamma$ and such that there exists some irreducible curve in $X$
  having numerical class in $\Q_{\,\geq
  0}\,\gamma$. 
We say that $\gamma$ is \emph{contractible} if 
there exist a toric variety $X_{\gamma}$ and
an equivariant morphism $\ph_{\gamma}\colon X\rightarrow  
X_{\gamma}$, with connected fibers, such that for every 
irreducible curve  $C$ in $X$,
$\ph_{\gamma}(C)=\{pt\}$ if and only if $[C]\in\Q_{\,\geq 0}\gamma$. 

\smallskip

This definition does not need the projectivity of $X$. 
We give a combinatorial characterization of contractibility in terms
of the fan of $X$, and we show that 
\emph{a class $\gamma$ is contractible if and
only if every irreducible invariant curve in the class is extremal in
every irreducible invariant
surface containing it} (theorem~\ref{guerra}).
In the projective case, this property is false for extremal classes:
  it can  
happen that every invariant curve in a class is extremal in every
  invariant  
subvariety containing it, but the class is not extremal in $X$ (see
example on  page~\pageref{esempio}).

When $X$ is projective,
all extremal classes are contractible, and
a contractible class $\gamma$ is extremal if and only if
$X_{\gamma}$ is projective. Hence, 
\emph{contractible non-extremal classes 
correspond to birational contractions to non-projective toric varieties}
(corollary~\ref{picasso}).
Moreover, we show that
a class $\gamma\in\NE(X)$ is
contractible if and only if it is extremal in 
the subvariety $A$
given by the intersection of all irreducible invariant divisors having
negative intersection with $\gamma$.

As mentioned above,
the main result of the paper is that 
\emph{when $X$ is projective, 
contractible classes span $\A(X)\cap\NE(X)$ as a semigroup} 
(theorem~\ref{bigteo}),  
namely every class in $\A(X)\cap\NE(X)$ decomposes as a linear
combination with  positive integral
coefficients of contractible classes. 
In the non-projective case, the situation is very different: 
if $\mathcal{C}$ is the set of contractible classes, in general
$\sum_{\gamma\in\mathcal{C}}\Q_{\,\geq 0}\,\gamma\subsetneqq\NE(X)$
(see remark on page~\pageref{rem}).

As an application of theorem~\ref{bigteo}, we show that 
\emph{when $X$ is projective, every curve having minimal degree with respect 
to some ample line bundle is extremal} (proposition~\ref{tie}).

This paper is about toric varieties. In the non-toric case, 
the situation is much more complicated; an account can be found in 
J.~Koll\'ar's survey paper~\cite{kollarff}. Briefly,
 if $X$ is a smooth projective
variety and $f\colon X\rightarrow Y$ a morphism with connected
fibers, there are essentially two reasons for which $Y$ can be
non-projective: either $f$ contracts a subcone of $\NE(X)$ 
which is not a face, 
or $f$ contracts a proper subset of some numerical class. In both
cases, this gives an effective 1-cycle in $Y$ which is homologous to
zero. In the toric case, we show that at least when $f$ is elementary, namely 
$\rho_X-\rho_Y=1$, the second case cannot happen 
(lemma~\ref{koll}),
essentially because numerical equivalence on $X$ implies numerical
equivalence on invariant subvarieties (see the remark on
page~\pageref{pin}). In~\cite{kollarff}, section~4, J.~Koll\'ar introduces
the notion of ``seemingly extremal ray'' in a smooth proper
algebraic space, in order to generalize the notion of (Mori) extremal ray.
Lemma~\ref{koll} easily implies that in a smooth complete toric
variety, contractible classes defined in this paper 
coincide with seemingly extremal rays that are contractible as
defined in~\cite{kollarff}.

The basic tools that we use are the language of primitive
collections and primitive relations, introduced by V.~V.~Batyrev 
(\cite{bat1,bat2}, see also \cite{sato}), and toric Mori theory
(M.~Reid \cite{reid}). Actually the combinatorial characterization of 
contractibility, 
in terms of the geometry of the fan, is 
implicitly present  already in \cite{reid}.

The structure of the paper is as follows: in
section 1  we briefly recall the definition and 
properties of primitive collections, and their link with toric Mori
theory. In section 2 we  define
contractible classes giving three equivalent conditions.
In section 3 we compare contractible and extremal classes when
$X$ is projective.
In section 4
we show that when $X$ is projective, 
contractible classes span $\A(X)\cap\NE(X)$ as a semigroup. 
Finally, in section 5, we study how contractible classes vary
under blow-up and blow-down, give some examples and show that if $X$ 
becomes projective after a single smooth equivariant blow-up,
then primitive relations span $\A(X)\cap\NE(X)$ as a semigroup. 

\medskip

\noindent\emph{Acknowledgments}.
This paper is part of my Ph.D.\ thesis. I am very grateful to my
advisor, Lucia Caporaso, and to Laurent Bonavero, for helping me
throughout the preparation of this work. I also wish to thank 
St\'ephane Druel for a careful reading of a preliminary version.

\stepcounter{subsubsec}
\subsubsection{Preliminaries}
For all the standard results in toric geometry, we refer to the books
of W.~Fulton~\cite{fulton} and T.~Oda~\cite{oda}.

Let $X$ be an $n$-dimensional toric variety: $X$ is described by a 
finite fan $\fx$ in the
vector space $N_{\Q}=N\otimes_{\Z}\Q$, where $N$ is a free abelian 
group of rank $n$.

We'll always assume $X$ smooth and complete, hence the support of
$\fx$ is the whole space $N_{\Q}$ and every cone in $\fx$ is generated
by a  part of a basis of $N$. However some of the results that we cite 
here hold more generally for a $\Q$-factorial complete toric variety,
as proposition~\ref{bacio} and Reid's results on toric Mori theory 
\cite{reid}.

If $x_1,\ldo,x_r\in N$, we will denote by $\langle
x_1,\ldo,x_r\rangle$ the cone in $N_{\Q}$ generated by $x_1,\ldo,x_r$,
namely:
\[ \langle
x_1,\ldo,x_r\rangle=\{\,
\sum_{i=1}^r\lambda_ix_i\,|\,\lambda_i\in\Q_{\,\geq 0}
\text{ for all }i=1,\ldo,r\,\}.\]
If $\sigma\in\fx$, we will always choose as a set of generators for
$\sigma$ a part of a basis of $N$.

We remember that for each $r=0,\ldo,n$ there is a bijection between
the cones of dimension $r$ in $\fx$ and the orbits of codimension $r$
in $X$; we'll denote by $V(\sigma)$ the closure of the orbit
corresponding to $\sigma\in\fx$ and $V(x)=V(\langle x\rangle)$ in case
of 1-dimensional cones. In what follows, we will 
refer to  the
subvarieties $V(\sigma)$ as invariant subvarieties. 

For each 1-dimensional cone $\rho\in\fx$, let $v_{\rho}\in\rho\cap N$
be its primitive generator, and 
\[  G(\fx)=\{v_{\rho}|\rho\in\fx\}  \]
the set of all generators in $\fx$.

\begin{defi}[V.~V.~Batyrev \cite{bat1}]
A subset $\lis{h}\subseteq G(\fx)$ is a \emph{primitive collection}
for $\fx$ if $\langle x_1,\ldo,x_h\rangle
\notin\fx$, but $\langle x_1,\ld,\check{x}_i,\ld,x_h\rangle\in\fx$
for each $i=1,\ldo,h$.
\end{defi}
\noindent We denote by $\PC(\fx)$ the set of all primitive collections
for $\fx$.   

Let $\sigma\in\fx$. We denote by $\RelInt\sigma$ the relative
interior of $\sigma$, namely the interior of $\sigma$ in its linear
span in $N_{\Q}$.
\begin{defi}
Let $P=\lis{h}\subseteq G(\fx)$ be a primitive collection. Since $X$ is 
complete, the point $\som{h}$ is contained in some cone of $\fx$; let 
$\sigma_P=\langle y_1,\ldo,y_k \rangle$ be the unique cone in $\fx$ such that
\[\som{h}\in \RelInt \sigma_P.
\]
Then we get a linear relation
\[
\som{h}-(a_1y_1+\cdots+ a_ky_k)=0
\]
with $a_i$ a positive integer for each $i=1,\ldo,k$. We call this
relation the \emph{primitive relation} associated to $P$.

The \emph{degree} of $P$ is the integer $\deg P = h-a_1-\cdots-a_k$.
\end{defi}

The set of the  primitive collections of a fan completely describes 
 the fan: an $n$-uple $\lis{n}\subseteq
G(\fx)$ generates a cone in the fan if and only if it doesn't
contain a primitive collection. 

Let $\A(X)$ be the group of algebraic 1-cycles on $X$ modulo numerical 
equivalence and $\N(X)=\A(X)\otimes_{\Z}\Q$.
We recall the well-known result:
\begin{prop}
\label{bacio}
The group $\A(X)$ is canonically isomorphic to the lattice of
integral relations among the elements of $G(\fx)$. A relation
\[ \sum_{x\in G(\fx)} a_x x=0, \qquad a_x\in\Z, \]
corresponds to a 1-cycle that has intersection $a_x$ with $V(x)$ for
all $x\in G(\fx)$.
\end{prop}
Hence, for
every primitive collection $P\in\PC(\fx)$, the associated primitive
relation defines a class $r(P)\in\A(X)$. 
Since the canonical class on $X$ is given by $K_X=-\sum_{x\in G(\fx)}
 V(x)$, for every primitive collection $P$ we have 
\[ -K_X\cdot r(P)=\deg P. \]

\noindent \emph{Notation:} we will often write primitive relations as
\[\som{h}=a_1y_1+\cdots+ a_ky_k \]
instead of $\som{h}-(a_1y_1+\cdots+ a_ky_k)=0$, i.\ e.\ writing
elements with negative coefficient on the right side. This must not 
be confused with the relation $-(\som{h})+ a_1y_1+\cdots+ a_ky_k=0$,
which is the opposite element in $\A(X)$.

\smallskip

\noindent \emph{Terminology:} by ``curve'' we always mean an 
\emph{irreducible} and \emph{reduced} curve in $X$.

\smallskip

Let  $\NE(X)\subset \N(X)$ the
cone of Mori, generated by classes of effective curves.
It is a result of M.~Reid~\cite{reid} that this cone is closed and polyhedral, 
generated by classes of invariant curves. The following simple lemma assures 
that primitive relations actually lie in $\NE(X)$:
\begin{lemma}
\label{lemma}
Let $\gamma\in\A(X)$ given by the relation
\[ a_1x_1+\cdots+a_hx_h-(b_1y_1+\cdots+b_ky_k)=0
\]
with $a_i,b_j\in\Z_{>0}$ for each $i,j$. If $\langle
y_1,\ldo,y_k\rangle\in\fx$, then $\gamma\in\NE(X)$.
\end{lemma}
\begin{dimo}
We have to show that $\gamma\cdot D\geq 0$ for every invariant nef divisor $D$
on $X$. All such divisors have the form
\[ D_{\ph}=-\sum_{x\in G(X)}\ph(x)V(x),
\]
where $\ph\colon N_{\Q}\rightarrow \Q$ is an upper 
convex support function, i.\ e.\
it is linear on each
cone of $\fx$, $\ph(N)\subseteq\Z$ and $\ph(u+v)\geq\ph(u)+
\ph(v)$ for each $u,v\in N_{\Q}$. Therefore
\begin{eqnarray*} 
D_{\ph}\cdot\gamma&=&-(a_1\ph(x_1)+\cdots+a_h\ph(x_h))+
b_1\ph(y_1)+\cdots+b_k\ph(y_k)
\\
&=&-( \ph(a_1x_1)+\cdots+\ph(a_hx_h))+\ph(b_1y_1+\cdots+b_ky_k)
\\
&=&-( \ph(a_1x_1)+\cdots+\ph(a_hx_h))+\ph( a_1x_1+\cdots+a_hx_h)\geq 0.
\end{eqnarray*}
\end{dimo}
\noindent\emph{Remark:}
a primitive relation does not need to be the numerical 
class of an  invariant curve $C$, and
the numerical class of an
invariant curve $C$ does not need to be a primitive relation.
Here is an example in dimension 3:

\vspace{10pt}

\hspace{30pt}\psfig{figure=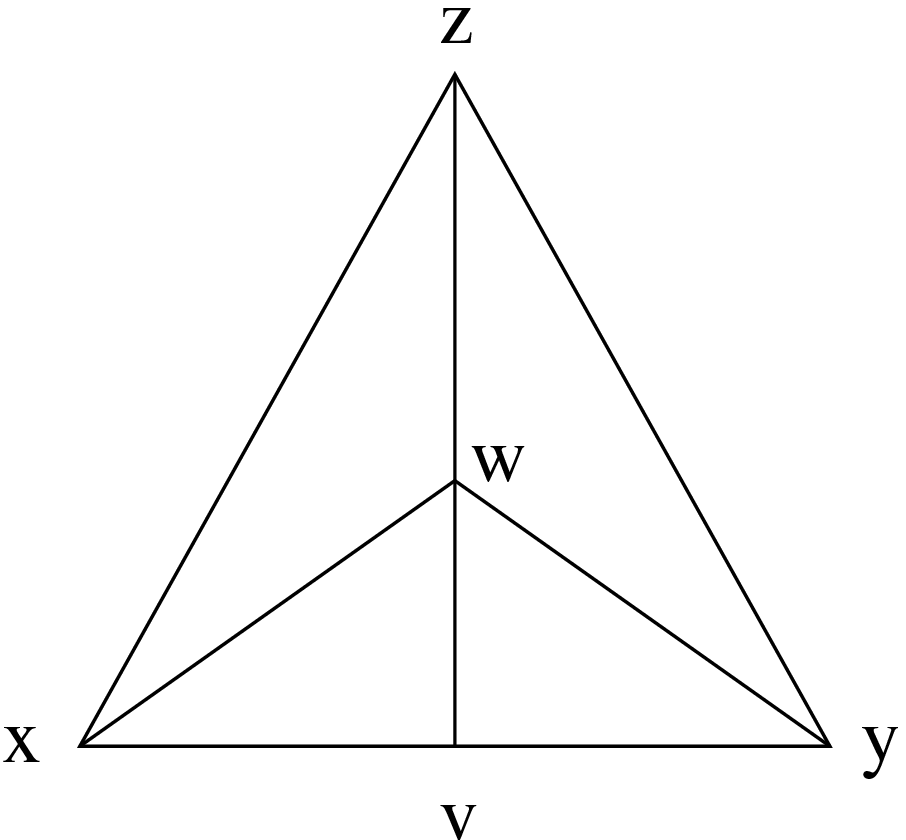,width=4cm}\hspace{40pt}\psfig{figure=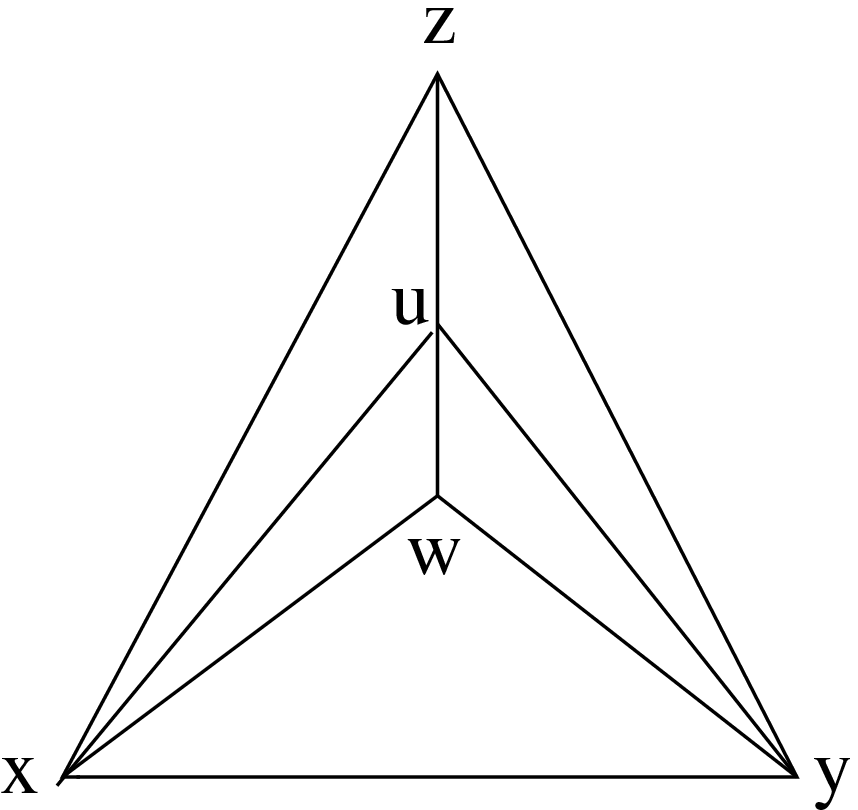,width=3.7cm} 

\vspace{10pt}
       
\noindent These are plane sections of two three-dimensional cones
$\langle x,y,z \rangle$; in every point corresponding to a ray we
indicate the generator of the ray. We suppose that $\{x,y,z\}$ is a
basis of the lattice and that $w=x+y+z$, $v=x+y$, $u=x+y+2z$.
On the left figure we have a curve $V(\langle z,w \rangle)$ whose relation
$x+y+z=w$ is not a primitive relation; on the right we have a
primitive collection $\{x,y,z\}$ whose primitive relation $x+y+z=w$ is
not associated to any invariant curve.

\smallskip

\noindent \emph{Remark: } let $\gamma$ be a primitive relation:
\[ \som{h}=a_1y_1+\cdots+a_ky_k. \]
Then, if there are invariant curves having numerical class $\gamma$,
they are necessarily contained in $V(\langle y_1,\ldo,y_k\rangle)
=V(y_1)\cap\cdots\cap V(y_k)$.

\smallskip

By Kleiman's criterion of ampleness~\cite{kleiman}, we know that $X$
is projective if and only if $\NE(X)$ is strictly convex.
In this case, its one-dimensional faces are called extremal rays.
We stress the fact that here extremal only refers to the geometry of
the cone $\NE(X)$; differently from Mori's extremal rays, we do not
require negative intersection with $K_X$.
 With the
following result, M.~Reid gives a precise description of the geometry of
the fan around a cone corresponding to a curve whose numerical class
lies in an extremal ray:
\begin{teo}[M.~Reid \cite{reid}, theorem 2.4] 
\label{reidextr}
Let $X$ be projective, $R$ an extremal ray of $\NE(X)$ and 
$\gamma\in R\cap\A(X)$ primitive along the ray. Then there exists a
primitive collection $P=\lis{h}$ such that $\gamma=r(P)$: 
\[\gamma:\quad \som{h}=a_1y_1+\cdots+a_ky_k.
\]
Moreover, for every $\nu=\langle z_1,\ldo,z_t\rangle$ 
such that $\{z_1,\ld,z_t\}\cap\{ x_1,\ld,x_h,y_1,\ld,y_k\}=\emptyset$
and $\langle y_1,\ld,y_k\rangle+\nu \in\fx$,
we have
\[\quad\langle
x_1,\ld,\check{x}_i,\ld,x_h,y_1,\ld,y_k\rangle+\nu
\in\fx\quad\text{for all }i=1,\ldo,h.
\]
\end{teo}
\noindent \emph{Remark:} suppose that $\gamma\in\NE(X)$ is either the
class of an invariant curve, or a primitive relation. Then, in both
cases, $\gamma$ is primitive in $\A(X)\cap\Q_{\,\geq 0}\gamma$. This is
because in both cases the relation has some coefficients equal to
1. In particular, when $R$ is an extremal ray of $\NE(X)$, there are a
unique class of an invariant curve and a unique primitive relation
contained in $R$, and they both coincide with the primitive element of
$R\cap\A(X)$. To avoid confusion, we will call extremal class (or
extremal curve) only such a primitive class.

\smallskip

As a consequence of theorem~\ref{reidextr}, 
we have an important description of
the cone of effective curves for projective toric varieties:
\begin{prop}[V.~V.~Batyrev \cite{bat1}]
Suppose that $X$ is projective. Then
the cone of effective curves $\NE(X)$ is generated by primitive relations.
\end{prop}

\noindent \emph{Remark:} if $X$ is non-projective, we just have  by
lemma~\ref{lemma}
\[ \sum_{P\in\PC(\fx)} \Q_{\,\geq 0}r(P)\subseteq\NE(X). \] 
We do not know if equality still holds.
Corollary~\ref{basta} tells us that equality holds if 
 $X$ becomes projective after a single smooth equivariant blow-up.

\stepcounter{subsubsec}
\subsubsection{Contractible classes}
We recall that $X$ is an $n$-dimensional, smooth, complete toric variety.

By a smooth equivariant blow-up, we mean the blow-up of a smooth toric
variety along a smooth, invariant subvariety. The resulting variety is
a smooth toric variety.

Let $f\colon X\rightarrow Y$ be a smooth equivariant blow-up along
$V(\tau)\subset Y$, with $\tau=\langle x_1,\ldo,x_h\rangle$.
Then $P=\lis{h}$ is a primitive collection in $\fx$, with relation
\[ r(P)\colon\quad \som{h}=x. \]
$V(x)$ is the exceptional divisor in $X$, and $r(P)$ is the class of a 
$\pr{1}$ contained in a fiber of $f$. 
We remark that $f$ contracts all 
irreducible curves in $X$ whose numerical
class is a multiple of $r(P)$.
In general, if $X$ is a non-toric smooth complex algebraic variety,
$f$ could contract only a proper subset of all the irreducible curves
having numerical class
in $\Q_{\,\geq 0}r(P)$; anyway, in such a case 
$Y$ would not be algebraic
 (see~\cite{bontak}~1.3, 
\cite{kollarff}~4.1.3, \cite{hartshorne} \S\ 3 in Appendix B
 and lemma~\ref{koll} in this paper).

The class $r(P)$ is not necessarily
extremal in $\NE(X)$:
\begin{prop}[L.~Bonavero \cite{bonavero}]
\label{bon}
Suppose that $X$ is projective. Then $r(P)$  is extremal in $\NE(X)$ 
if and only if $Y$ is projective.
\end{prop}
We remark that even if $r(P)$ is not extremal, it is ``contractible''
in the sense that there exists a morphism
$f\colon X\rightarrow Y$ such that for any irreducible curve $C\subset
X$, $f(C)$ is a point if and only if 
$[C]\in\Q_{\,\geq 0}r(P)$. 

We want to give an example of this behaviour. \label{esempio}
Let $Y$ be the only complete
non projective toric
3-fold with Picard number 4 (see \cite{oda}, page
85, and also \cite{bonavero}). 
The fan of $Y$ is on the left side of the figure (it is projected onto
the face $\langle e_1,e_2,e_3\rangle$). The set $\{e_1,e_2,e_3\}$ is a
basis of the lattice, $e_0=-e_1-e_2-e_3$ and $f_i=e_0+e_i$ for
$i=1,2,3$.

\vspace{10pt}
$\fy\quad
$\psfig{figure=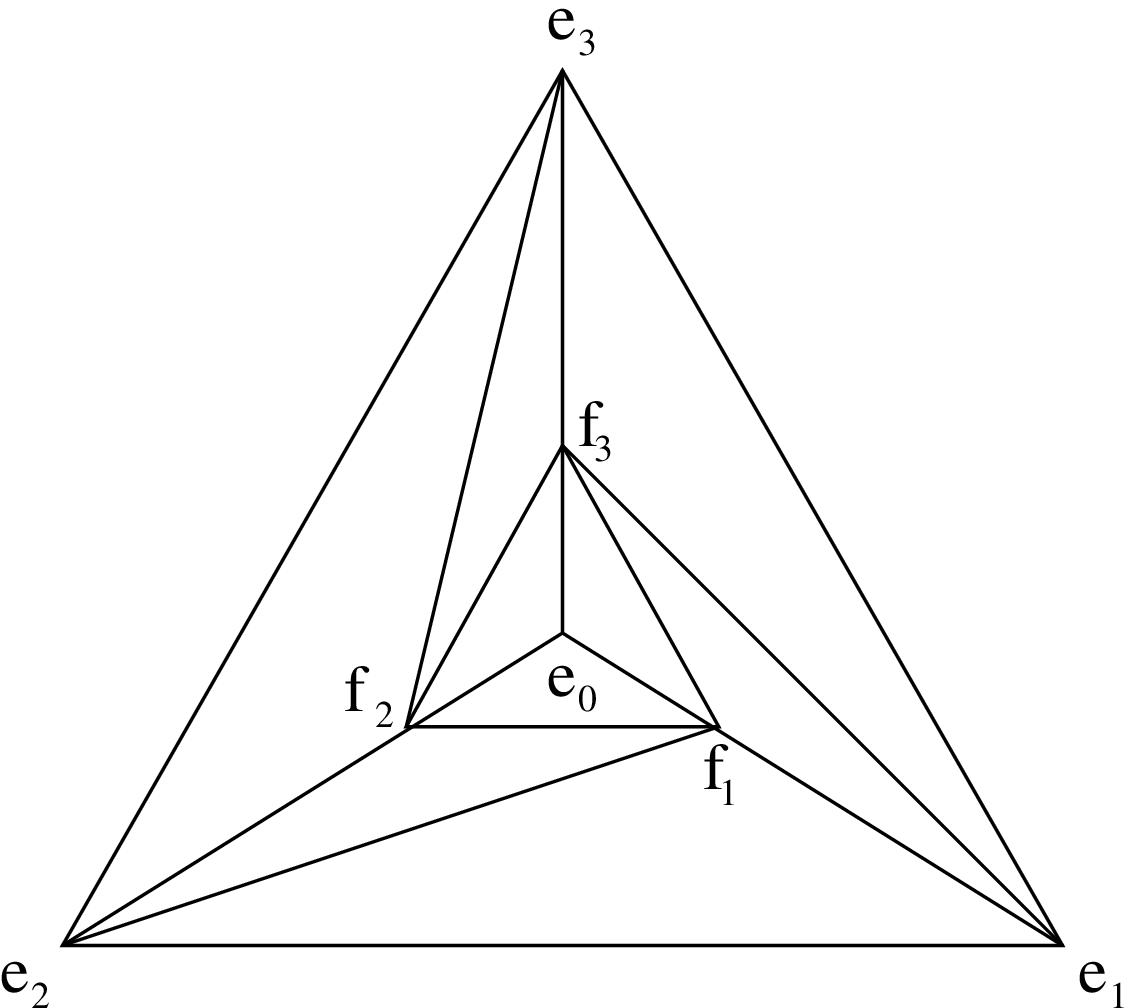,width=4.6cm}\hspace{30pt}$\fx\quad$\psfig{figure=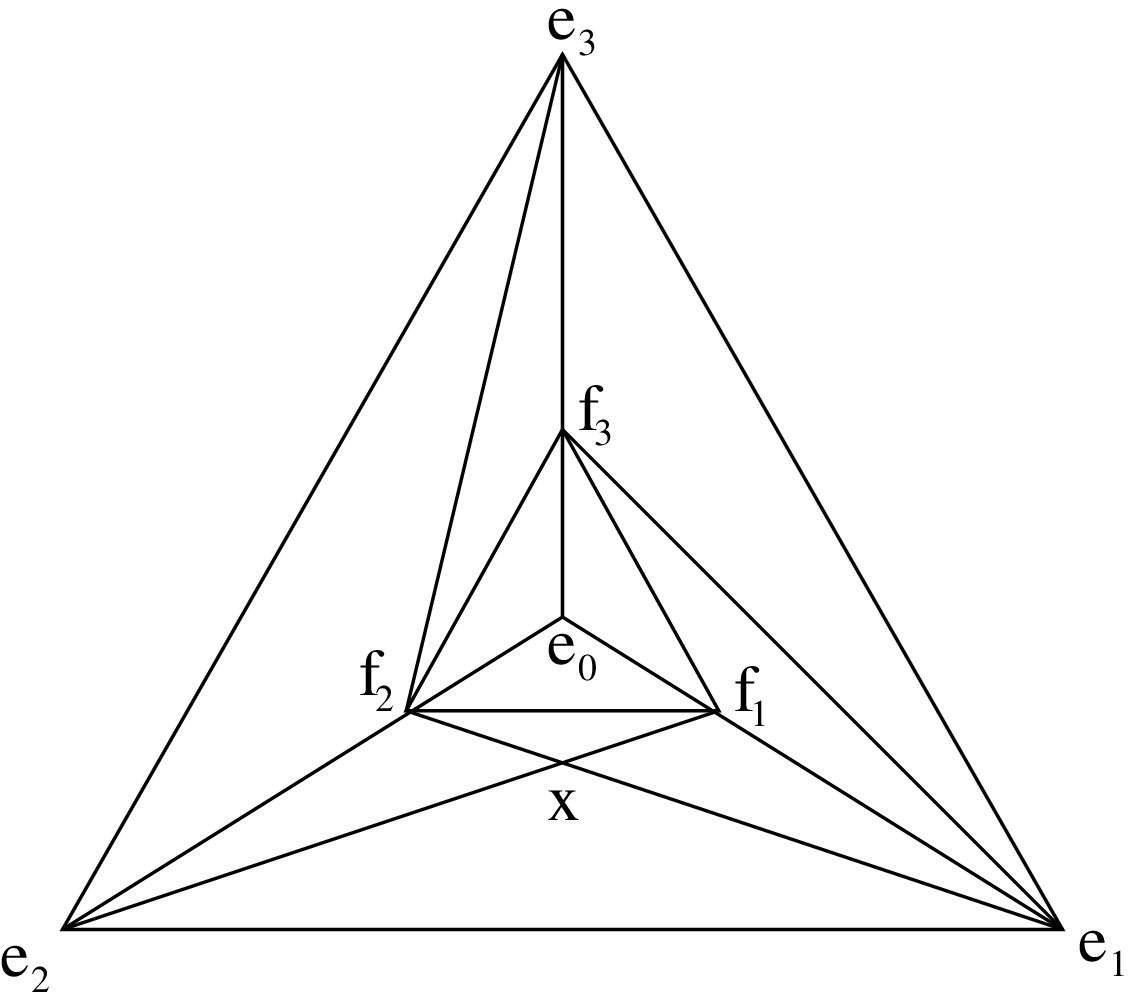,width=4.6cm} 
\vspace{10pt}

\noindent Let $X$ be the blow-up of $Y$ along the curve $V(\langle
f_1,e_2\rangle)$. The fan of $X$ is on the right side of the figure. 
The set $\{f_1,e_2\}$ becomes a primitive collection in $\fx$, with
relation $f_1+e_2=x$. Let $\gamma\in\NE(X)$ be the class
corresponding to this relation.
It is easy to see that $X$ is projective: it is obtained from $\pr{3}$
by 4 smooth equivariant blow-ups. By proposition~\ref{bon}, $\gamma$
is not extremal in $\NE(X)$, otherwise $Y$ would be projective. 
The variety $X$ has Picard number 5, and its extremal classes are given by the
primitive relations:
\[ \begin{array}{crcl}
\omega_1:&\qquad e_1+f_2&=&x\\
\omega_2:&\qquad f_1+e_3&=&e_1+f_3\\
\omega_3:&\qquad e_2+f_3&=&f_2+e_3\\
\omega_4:&\qquad e_0+x\,&=&f_1+f_2\\
\omega_5:&\qquad f_1+f_2+f_3&=&2e_0.
\end{array} \]
$\NE(X)$ is a simplicial cone in a 5-dimensional vector space,
and $\gamma=\omega_1+\omega_2+\omega_3$. The invariant curves in $X$
having class $\gamma$ are $C_1=V(\langle e_1,x\rangle)=V(e_1)\cap V(x)$
 and
$C_2=V(\langle f_2,x\rangle)=V(f_2)\cap V(x)$.  
Let's consider the irreducible invariant surfaces containing $C_1$ or $C_2$.
The surface 
$V(x)$ is $\pr{1}\times\pr{1}$ and
$C_1$, $C_2$ in $V(x)$ are the invariant fibers of one of the two
fibrations on $\pr{1}$. The surface 
$V(e_1)$ is $\pr{2}$ blown-up in two points,
and $C_1$ is one of the exceptional curves. The surface $V(f_2)$ is
obtained from $\pr{2}$ blown-up in two points, blowing-up a fixed
point lying on one exceptional curve $E$; in $V(f_2)$ the curve $C_2$ is the
strict transform of $E$ and it has self-intersection $-2$, hence it is
extremal.
Thus we see that even if $\gamma$ is not
extremal in $X$, $C_1$ and $C_2$ are extremal in both irreducible invariant
surfaces containing them.

We are going to define a class of elements in $\NE(X)$ with properties
similar to those of $\gamma$.
\begin{teo}
\label{guerra}
Let  $X$ be a smooth, complete toric variety and
let $\gamma\in\A(X)\cap\NE(X)$. Then the following conditions are equivalent:

\smallskip

\noindent (i) $\gamma$ is the class of some invariant curve in $X$,
and every invariant curve having class $\gamma$ is extremal in every
irreducible invariant surface containing it.

\smallskip

\noindent (ii) $\gamma$ is a primitive relation of the form 
\[ \som{h}=a_1y_1+\cdots+a_ky_k \]
and for every $\nu=\langle z_1,\ld,z_t\rangle\in\fx$  such that 
$\{z_1,\ld,z_t\}\cap\{x_1,\ld,x_h,y_1,\ld,y_k\}=\emptyset$ and
$\langle y_1,\ldo,y_k\rangle +\nu\in\fx$, 
 then
\[ \langle x_1,\ld,\check{x}_i,\ld,x_h,y_1,\ld,y_k\rangle+\nu\in\fx
\quad\text{ for all }i=1,\ldo,h.\]

\smallskip

\noindent (iii)  $\gamma$ is primitive in $\A(X)\cap\Q_{\,\geq
  0}\gamma$, there exists some irreducible curve having numerical
  class in $\Q_{\,\geq 0}\gamma$,  and
there exist a
toric variety $X_{\gamma}$ and an equivariant morphism 
$\ph_{\gamma}\colon X
\rightarrow X_{\gamma}$, with  connected fibers,
such that for every irreducible curve $C\subset X$
\[
\ph_{\gamma}(C)=\{pt\}\Longleftrightarrow [C]\in\Q_{\,\geq 0}\gamma. \]
\end{teo}
\begin{defi}
\label{pioggia}
We say that $\gamma$ is \emph{contractible} if one of the equivalent
conditions of theorem~\ref{guerra} is satisfied.
\end{defi}

\noindent\emph{Remarks:} 

\noindent $\bullet\ \,$  when $X$ is projective,
by Reid's theorem~\ref{reidextr}, an extremal
class is always contractible.

\noindent $\bullet\ $ We recall that all complete smooth toric
surfaces are projective.
For surfaces, the
notions of contractible class and extremal class coincide. 

\noindent $\bullet\ $ As we saw in the preceeding example, if
$X$ is projective, a class can be extremal in all irreducible
invariant hypersurfaces containing some 
invariant curve of the class, without being
extremal.

\noindent $\bullet\ $ \label{pin}
If $V\subset X$ is an irreducible invariant
subvariety, then the natural map $\Pic X\rightarrow\Pic V$ induced by
restriction of divisors is surjective, hence  we have an inclusion 
$\N(V)\subseteq\N(X)$. In other words, if two curves
$C_1,C_2\subset V$ are numerically equivalent in $X$, then the same
holds in $V$. Thus, if $\ph_{\gamma}$ is a morphism as in $(iii)$,
then there are two possibilities:
\begin{enumerate}[a.]
\item $V$ does not contain any 
irreducible curve having numerical class $\gamma$,
  hence $\ph_{\gamma\,|V}\colon V\rightarrow \ph_{\gamma}(V)$ 
is an isomorphism;
\item $\gamma$ defines a class in $\NE(V)$ and  
$\ph_{\gamma\,|V}\colon V\rightarrow \ph_{\gamma}(V)$
is the contraction of this class.
\end{enumerate}
This is a fundamental property which does not hold for general (non-toric)
varieties. 

\smallskip

The proof of the two implications $(i)\Rightarrow (ii)$ and
$(ii)\Rightarrow (iii)$ of theorem~\ref{guerra} is based on Reid's
proof of theorem~2.4 in~\cite{reid}; looking carefully at Reid's
construction, it actually turns out that to get the existence of a
morphism as in $(iii)$, you just need that
invariant curves in the class are extremal in all invariant 
surfaces containing them.
We also remark that reformulations of Reid's theorem~\ref{reidextr} 
in terms of
primitive relations are already present in Batyrev's and 
Sato's work
(see \cite{bat2}, section 2.3, and \cite{sato}, sections 4 and 5).
\begin{dimo}[Proof of theorem \ref{guerra}]

\noindent $(i)\Rightarrow (ii)\ $
Let $C\subset X$ be an invariant curve having numerical class $\gamma$ and let
\[ \langle 
e_1,\ldo,e_{n-1}\rangle=\langle e_1,\ldo,e_{n-1},e_n \rangle\cap \langle
e_1,\ldo,e_{n-1},e_{n+1}\rangle \]
be the associated cone. Then the relation associated to $[C]$ is
\[ \sum_{i=1}^{n-1}a_ie_i +e_n+e_{n+1}=0, \qquad a_i\in\Z. \]
For each $i=1,\ldo,n-1$ we consider the invariant surface $S_i$
corresponding to $\langle e_1,\ld,\check{e}_i,\ld,e_{n-1}\rangle$.
Let $N_i$ be the subgroup of $N$ generated by $e_1,\ld,\check{e}_i, 
\ld,e_{n-1}$; we denote by $\overline{x}$ 
the image  of an element $x\in N$
in the quotient group $N/N_i$.
In $\f_{S_i}$ the curve $C$ corresponds to the ray $\langle
\overline{e}_i\rangle$ and has relation 
$a_i \overline{e}_i +\overline{e}_n+\overline{e}_{n+1}=0$. 
Since $C$ is extremal in $\NE(S_i)$, we have:

\noindent $\bullet\ $
 if $i$ is such that $a_i>0$, then 
$a_i=1$ and $S_i\simeq\pr{2}$, so $\langle
e_1,\ld,\check{e}_i,\ld,e_{n+1}\rangle\in\fx$; 

\noindent $\bullet\ $ if $i$ is such that $a_i=0$,
then $S_i$ is a Hirzebruch surface $\mathbb{F}_a$, so
there exists $e_i'\in G(\fx)$ such that
$\langle e_1,\ld,\check{e}_i,e_i',\ld,e_{n-1},e_n\rangle$,
$\langle
e_1,\ld,\check{e}_i,e_i',\ld,e_{n-1},e_{n+1}\rangle$ are in
$\fx$.

We reorder $e_1,\ldo,e_{n-1}$ in such a way that $a_1,\ldo,a_q$
are negative and $a_p,\ldo,a_{n-1}$ are positive, with 
$0\leq q<p\leq n$.
Then the relation associated
to $[C]$ is 
\[e_p+\cdots +e_{n+1}=b_1e_1+\cdots +b_qe_q, \] 
with $b_i=|a_i|$ for $i=1,\ldo,q$; moreover, it is a primitive relation.

We remark that if $p=q+1$, then $(ii)$ is already proved; thus let's
assume in the sequel $p>q+1$.

Let $
\Delta_{\gamma}\subset\fx$ be the set of cones
$\nu=\langle z_{q+1},\ld,z_{p-1}\rangle$ such that the cone  
$\langle
e_1,\ld,e_q,z_{q+1},\ld,z_{p-1},e_p,\ld,\check{e}_i,\ld,e_{n+1}\rangle$
has dimension $n$ and is in $\fx$ for all $i=p,\ld,n+1$.
We have just proved that
$\langle e_{q+1},\ldo,e_{p-1}\rangle\in \Delta_{\gamma}$, thus 
$\Delta_{\gamma}$ is
non-empty. 
For each $\nu=\langle z_{q+1},\ld,z_{p-1}\rangle\in\Delta_{\gamma}$,
we get invariant curves having class $\gamma$, corresponding to the cones 
\[ \langle
e_1,\ld,e_q,z_{q+1},\ld,z_{p-1},e_p,\ld,\check{e}_i,\ld,\check{e}_j,\ld,
e_{n+1}\rangle, \qquad i,j\in\{p,\ld,n+1\}.\]
Again, looking at the irreducible invariant surfaces containing these curves,
we get that for all $l=q+1,\ld,p-1$ there exists
$z_l'\in G(\fx)$ such that $\nu_l'=\langle
z_{q+1},\ld,\check{z}_l,z_l',\ld, z_{p-1}\rangle$ is in $\Delta_{\gamma}$.

 Let $U=\langle e_1,\ld,e_q,e_p,\ld,e_{n+1}\rangle$ and
$P(\nu)=U+\nu$ for all $\nu\in\Delta_{\gamma}$.
Each $P(\nu)$ is a union of $n$-dimensional cones of
$\fx$. There are three types of
$(n-1)$-dimensional cones of $\fx$ lying
in $P(\nu)$:

\noindent $\bullet\ $ the cones of type $\langle
  e_1,\ld,\check{e}_j,\ld,e_q,e_p,\ld,\check{e}_i,\ld,e_{n+1}\rangle+\nu $, 
with $j\in\{1,\ld,q\}$ and $i\in\{p,\ld,n+1\}$,
which are external faces of $P(\nu)$; 

\noindent $\bullet\ $ the cones of type  $\langle
  e_1,\ld,e_q,e_p,\ld,\check{e}_i,\ld,\check{e}_j,\ld,e_{n+1}\rangle+\nu $,
with $i,j\in\{p,\ld,n+1\}$, which are internal to $P(\nu)$, because
they are intersection of two $n$-dimensional cones in $P(\nu)$;

\noindent $\bullet\ $ the cones of type 
 $\langle
  e_1,\ld,e_q,z_{q+1},\ld,\check{z}_j,\ld,z_{p-1},e_p,\ld,\check{e}_i,
  \ld,e_{n+1}\rangle$ with
$j\in\{q+1,\ld,p-1\}$ and $i\in\{p,\ld,n+1\}$; these cones 
are intersection of two $n$-dimensional
  cones, one  in $P(\nu)$ and one in $P(\nu_j')$, where $\nu_j'=\langle
z_{q+1},\ld,\check{z}_j,z_j',\ld, z_{p-1}\rangle$.

Consider now  the union $\mathcal{V}$
of all $P(\nu)$
for $\nu\in\Delta_{\gamma}$. We have just shown that  $\mathcal{V}$ is
a union of $n$-dimensional cones of $\fx$ and that its external faces do
not contain $\langle e_1,\ldo,e_q\rangle$. Hence, if $\langle
z_{q+1},\ld,z_{p-1}\rangle$ is such that $\{z_{q+1},\ld,z_{p-1}\}\cap
\{e_1,\ld,e_q,e_p,\ld,e_{n+1}\}=\emptyset$ and $\langle e_1,\ld,e_q,
z_{q+1},\ld,z_{p-1}\rangle\in\fx$, it must be  $\langle
z_{q+1},\ld,z_{p-1}\rangle\in \Delta_{\gamma}$.
Hence we have shown that $(ii)$ holds for $\gamma$.

\medskip

\noindent $(ii)\Rightarrow (iii)\ $ $\gamma$ is a primitive relation of
the form
\[ \som{h}=a_1y_1+\cdots+a_ky_k. \]
For $\nu=\{0\}$, we get $\langle
x_1,\ld,\check{x}_i,\ld,x_h,y_1,\ld,y_k\rangle\in\fx$ for all
$i=1,\ldo,h$.

Let $\Delta_{\gamma}\subset\fx$ be the set of cones $\nu=\langle
z_1,\ld,z_t\rangle$  of
dimension $t=n-h-k+1$ such that
$\{z_1,\ld,z_t\}\cap\{x_1,\ld,x_h,y_1,\ld,y_k\} =\emptyset$ and
$\nu+\langle y_1,\ld,y_k\rangle\in\fx$. 
Any $n$-dimensional cone of $\fx$ containing $\langle
y_1,\ldo,y_k\rangle$ must have the form $\langle x_1,\ld,\check{x}_i,
\ld,x_h,y_1,\ld,y_h,z_1,\ld,z_t\rangle$ for some $i\in\{1,\ld,h\}$,
hence $\Delta_{\gamma}$ is non-empty.

Let $\nu=\langle z_1,\ld,z_t\rangle\in\Delta_{\gamma}$ 
and fix $i\in\{1,\ld,h\}$, $j\in\{1,\ld,t\}$. 
The cone $\langle x_1,\ld,\check{x}_i,\ld,x_h,y_1,
\ld,y_k,z_1,\ld,\check{z}_j, \ld,z_t\rangle$ has
dimension $n-1$ and must be the intersection of two $n$-dimensional
cones in $\fx$; therefore there exists a $z_j'\in G(\fx)$ such that
$\nu_j'=\langle
z_1,\ld,\check{z}_j,z_j',\ld,z_t\rangle\in\Delta_{\gamma}$.  

Set $U=\langle x_1,\ld,x_h,y_1,\ld,y_k\rangle$, 
$P(\nu)=U+\nu$ for $\nu\in\Delta_{\gamma}$ and let  $\mathcal{V}$ be 
the union of all $P(\nu)$. The same description given in the first
part of the proof holds: distinct $P(\nu)$
intersect each other along common faces; $\mathcal{V}$ is a union of
$n$-dimensional cones of $\fx$ and contains all $(n-1)$-dimensional 
cones corresponding to invariant curves having numerical class $\gamma$.
 
Suppose that $k>0$, i.\ e.\ that $\gamma$ is not numerically
effective. 
Then the boundary of $\mathcal{V}$ is a union of faces in $\fx$, none
of which contains $\langle y_1,\ldo,y_k\rangle$. 
The set
\[ \f_{X_{\gamma}}=\fx\smallsetminus\{\tau\in\fx\,|\,\dim\tau=n-1,\
V(\tau)\text{ has numerical class }\gamma\} \] 
is a fan, different from $\fx$
only inside $\mathcal{V}$; inside $\mathcal{V}$,
its $n$-dimensional cones are the $P(\nu)$, for $\nu\in\Delta_{\gamma}$.  
The fiber over the point $V(P(\nu))\in X_{\gamma}$ is $V(\nu+\langle
y_1,\ldo,y_k\rangle)\simeq\pr{h-1}$.

If $k=0$, i.\ e.\ the class $\gamma$ is numerically effective, then
$U$ is a linear subspace of $N_{\Q}$, and $\mathcal{V}$ is
a polyhedral decomposition of $N_{\Q}$.
If you remove from $\fx$ all $(n-1)$-dimensional cones corresponding to
curves having class $\gamma$, you get a degenerate fan with vertex 
$U$; then taking the quotient of each cone by $U$ you get the fan of
$X_{\gamma}$. 

\smallskip

\noindent $(iii)\Rightarrow (i)\ $ First, we show that the
implication is true when $\dim X=2$. 
In this case $X$ is obtained by a finite sequence of blow-ups at fixed
points from $\pr{2}$ or from a Hirzebruch surface $\mathbb{F}_a$ (see
\cite{oda}, page 42); in particular, $X$ is projective.
There are only three possibility for $\ph_{\gamma}$:

\noindent $\bullet\ $ $\ph_{\gamma}$ contracts $X$ to a point; then
$X\simeq\pr{2}$ and $\gamma$ is the class of a line;

\noindent $\bullet\ $ $\ph_{\gamma}\colon X\rightarrow\pr{1}$; then
$X$ is a Hirzebruch surface $\mathbb{F}_a$ and $\gamma$ is the class
of a fiber;

\noindent $\bullet\ $ $\ph_{\gamma}$ is birational and its
exceptional locus is a curve $E$ with negative self-intersection; $E$
is the only curve having class $\gamma$.

In all three cases, we see that $\gamma$ is extremal in $\NE(X)$.

When $\dim X>2$, consider an irreducible invariant surface $S\subset
X$. If $S$ contains some curves having class $\gamma$, then 
as remarked on page \pageref{pin}, the morphism 
$\ph_{\gamma\,|S}\colon S\rightarrow \ph_{\gamma}(S)$ is the
contraction of a contractible class in $\NE(S)$; hence
this class is extremal in $\NE(S)$ and the statement follows.
\end{dimo}

\begin{cor}
\label{pluto}
Let $X$ be a complete, smooth toric variety,
$\gamma\in\NE(X)$ a contractible class and 
$\ph_{\gamma}\colon X
\rightarrow X_{\gamma}$ the associated morphism.

Suppose that $\gamma$ is numerically effective:
\[\som{h}=0. \]
Then $X_{\gamma}$ is smooth of dimension $n-h+1$ and
$\ph_{\gamma}$ is a $\pr{h-1}$-bundle.

Suppose that $\gamma$ is not numerically effective: 
\[\som{h}=a_1y_1+\cdots+a_k y_k, \quad k>0. \]
Then $\ph_{\gamma}$ is birational, with
  exceptional loci $A\subset X$, $B\subset
X_{\gamma}$ given by $A=V(\langle
y_1,\ldo,y_k\rangle)$, $B=V(\langle x_1,\ld,x_h,
y_1,\ld,y_k\rangle)$; $\dim A= n-k$, $\dim B= n-h-k+1$ and 
$\ph_{\gamma|A}\colon A\rightarrow B$ is a $\pr{h-1}$-bundle.
Moreover, $X_{\gamma}$ is simplicial if and only if $k=1$; it is
smooth if and only if $k=1$ and $a_1=1$, in which case $\ph_{\gamma}$
is a smooth equivariant blow-up.
\end{cor}
\noindent\emph{Remark:} here with ``projective bundle'' we mean 
locally trivial
on Zariski open subsets; hence it is actually the projectivized of some
vector bundle on the base (see~\cite{hartshorne}, 7.10, page 170).

\begin{proof}
In the case of $\gamma$ numerically effective, the fact that 
$\ph_{\gamma}$ is a $\pr{h-1}$-bundle follows from the combinatorial
description in theorem~\ref{guerra}, $(ii)$, by standard results of
toric geometry (see~\cite{ewald2}, theorem 6.7 on page 246,
or~\cite{fulton}, page 41).

When $\gamma$ is not numerically effective, the only unclear point 
in the statement is that $\ph_{\gamma|A}\colon A\rightarrow B$ is a
$\pr{h-1}$-bundle. But this follows from the preceeding case, because
in $A$ we are contracting a numerically effective class.
\end{proof}

The following corollary points out that it is enough to check
contractibility in a suitable subvariety:
\begin{cor}
\label{paris}
Let $X$ be a complete, smooth toric variety and
$\gamma\in\A(X)\cap\NE(X)$.
Let $D_1,\ldo,D_k$ be the irreducible invariant divisors in $X$ having
negative intersection with $\gamma$ and let $A=D_1\cap\cdots\cap D_k$.
Then $\gamma$ 
is contractible in $X$ if and only if it is contractible in
$A$.
\end{cor}
\begin{dimo}
All irreducible curves in $\gamma$ are contained in $A$ and they are still
numerically equivalent in $A$ (see the fourth remark after 
theorem~\ref{guerra}, on page~\pageref{pin}).
Hence $\gamma$ actually defines a class in $\NE(A)$.

Consider the relation associated to $\gamma$ in $X$:
\[ b_1x_1+\cdots+b_hx_h-(a_1y_1+\cdots+a_ky_k)=0\ \text{ with
  }b_i,a_j> 0\text{ for all }i,j.
\]
Then $\{D_1,\ldo,D_k\}=\{V(y_1),\ldo,V(y_k)\}$ and 
$A=V(\langle y_1,\ldo,y_k\rangle)$.
The primitive relation associated to $\gamma$ in $A$ is 
\[ b_1\overline{x}_1+\cdots+b_h\overline{x}_h=0,\]
where $\overline{x}_i$ is the image of $x_i$ in the quotient group of
$N$ by the subgroup spanned by $y_1,\ldo,y_k$. 
Condition $(ii)$ of theorem~\ref{guerra} for
$\gamma$ gives
the same in $X$ and in $A$, hence the statement follows.
\end{dimo}

The following lemma shows that an equivariant morphism $f\colon
X\rightarrow Y$ with connected fibers and such that $\rho_X-\rho_Y=1$,
is always of type $\ph_{\gamma}$ for a contractible class $\gamma$.
\begin{lemma}
\label{koll}
Let $X$ be a complete, smooth toric variety, $Y$ a toric variety and
$f\colon X\rightarrow Y$ an equivariant morphism with connected
fibers. Suppose that there exists $\gamma\in\NE(X)\cap\A(X)$,
primitive in $\A(X)\cap\Q_{\,\geq 0}\gamma$, such that for every 
irreducible curve
$C\subset X$ such that $f(C)=\{pt\}$, we have $[C]\in\Q_{\,\geq
  0}\gamma$. Then $\gamma$ is contractible and $f=\ph_{\gamma}$.
\end{lemma}
\begin{dimo}
We have to show that if $C\subset X$ is 
an irreducible curve such that $[C]\in\Q_{\,\geq
  0}\gamma$, then $f(C)=\{pt\}$.

Let $D_1,\ldo,D_k$ be the irreducible
invariant divisors in $X$ having negative intersection with $\gamma$ 
and let $A=D_1\cap\cdots\cap D_k$. Considering the restriction
$f_{|A}\colon A\rightarrow f(A)$, we can suppose that $\gamma$ is
numerically effective and $A=X$.

Moreover, it is enough to prove the statement for the restriction
$f_{|S}\colon S\rightarrow f(S)$ to any irreducible invariant surface
$S\subset X$,  thus we can suppose $\dim X=2$.

Since every curve contracted by $f$ is numerically effective, 
it must be $\dim Y<\dim X$ or $f$ isomorphism. If $\dim Y=0$, the
statement is clear. If $\dim Y=1$, we have a fibration $f\colon
X\rightarrow\pr{1}$, and $\gamma$ is the class of the generic fiber.
 Since $f$ cannot contract any exceptional curve, it must
 be a toric $\pr{1}$-bundle: hence $\gamma$ is extremal in $\NE(X)$ 
and $f=\ph_{\gamma}$.
\end{dimo}
As we already remarked in the introduction and 
at the beginning of section~2, in the non-toric
case the statement of lemma~\ref{koll} is false if $Y$ is
non-projective; in fact, it may happen that $f\colon X\rightarrow Y$
contracts only a proper subset of the numerical class~$\gamma$.

\stepcounter{subsubsec}
\subsubsection{Contractible versus extremal in projective varieties}
In this section we characterize contractible, non-extremal classes in
a projective toric
variety. Moreover, we give a combinatorial criterion for
a primitive relation to be contractible.

For numerically effective classes, the 
following result is immediate  from
corollary~\ref{pluto}:
\begin{cor}
Suppose that $X$ is projective and $\gamma\in\NE(X)$  
is contractible and numerically effective.
Then $\gamma$ is extremal and $X_{\gamma}$ is projective.
\end{cor}
Hence, by corollary~\ref{paris}, when $X$ is projective
a class $\gamma\in\NE(X)$ is
contractible if and only if it is extremal in the subvariety $A$
given by the intersection of all irreducible invariant divisors having
negative intersection with $\gamma$.

For non numerically effective contractible classes we have 
the following:
\begin{lemma}
Let $\gamma$ be a contractible class, non numerically effective,
and let $\ph_{\gamma}\colon X\rightarrow X_{\gamma}$ be the associated
birational morphism.

Then $X_{\gamma}$ is projective if and only if $X$ is projective and
$\gamma$ is an extremal class.
\end{lemma}
\begin{dimo}
Suppose $X_{\gamma}$ that is projective and let $H\in\Pic
X_{\gamma}$ ample. Then for every irreducible
curve $C\subset X$ we have
$\ph_{\gamma}^{*}(H)\cdot C\geq 0$, and
\[\ph_{\gamma}^{*}(H)\cdot C=0 \ \Leftrightarrow\ 
\ph_{\gamma}(C)=\{pt\}\ \Leftrightarrow\ [C]\in\Q_{\,\geq 0}\gamma. \]
Hence $\gamma$ generates an extremal ray in $\NE(X)$, which also
implies that $\NE(X)$ is strictly convex and $X$ is projective. 

Viceversa,  if $X$ is projective and $\gamma$ is extremal, it is known
that $X_{\gamma}$ is projective.
\end{dimo}

\begin{cor}
\label{picasso}
If $X$ is projective, then contractible non extremal classes in $\NE(X)$
correspond to birational contractions to non projective varieties.
\end{cor}

We end this section with a combinatorial
criterion for contractibility, which gives a simple combinatorial
algorythm to determine, given all primitive relations in $\fx$,
which are the contractible ones. It's remarkable that, when $X$ is
projective, there is no analogous algorythm to determine which 
primitive relations are  extremal in $\NE(X)$.

Proposition~\ref{criterio} has been proven 
by H.~Sato (\cite{sato}, theorem 4.10) for a
primitive relation $r(P)$ corresponding to
a smooth equivariant blow-down; the same proof holds for general
primitive relations.
\begin{prop}
\label{criterio}
Let $P=\{x_1,\ldo,x_h\}$ be a primitive collection in $\fx$, with
primitive relation $r(P):\ x_1+\cdots +x_h=a_1y_1+\cdots +a_ky_k$. Then
$r(P)$ is contractible if and only if for every primitive collection
$Q$ of $\fx$ such that
$Q\cap  P\neq\emptyset$ and $Q\neq P$, the set $(Q\smallsetminus 
P)\cup\{y_1,\ldo,y_k\}$ contains a primitive collection.
\end{prop}
In particular, when a primitive relation $r(P)$ is numerically
effective, it is contractible if and only if $P$ is disjoint from all 
other primitive collections of $\fx$.

\stepcounter{subsubsec}
\subsubsection{A property of $\NE(X)$ for $X$ projective}
In this section, we prove the 
\begin{teo}
\label{bigteo}
Let $X$ be a smooth projective toric variety. Then
for every $\eta\in\A(X)\cap\NE(X)$ there is a decomposition
\[ \eta=m_1 \gamma_1+\cdots+m_r \gamma_r \]
with $\gamma_i$ contractible and $m_i\in\Z_{>0}$ for all $i=1,\ldo,r$.
\end{teo}
We will prove theorem \ref{bigteo} by reducing the analysis to
invariant surfaces of $X$.
We remark that, since every curve on $X$ is numerically equivalent to
a linear 
combination of invariant curves with positive integral coefficients
(see \cite{reid}, proposition (1.6)),
it is enough to prove the theorem for classes of invariant curves.

In dimension 2, contractible classes coincide with extremal classes,
because all complete toric surfaces are projective.
Therefore, theorem \ref{bigteo} can be restated for surfaces as
follows:
\begin{prop}
\label{lucia}
Let $S$ be a smooth projective toric surface. Then every class 
$\eta \in\A(S)\cap\NE(S)$ can be decomposed as a linear combination with
positive integral coefficients of extremal classes.
\end{prop}
\begin{dimo}[Proof of proposition \ref{lucia}]
We have to show:\\
(P)\ \ \emph{if $C\subset S$ invariant, then $C\equiv\sum_{i}m_iC_i$, 
with $m_i\in\Z_{>0}$ and $C_i$ extremal \hspace*{6 mm} in $\NE(S)$ 
for all $i$.}

(P) is true for minimal toric surfaces, namely for $\pr{2}$ and for
the Hirzebruch surfaces $\mathbb{F}_a$ (see for instance~\cite{oda},
page 108). We show that if $\pi\colon
S\rightarrow T$ is an equivariant smooth blow-up, and (P) is true for
$T$, then it is true for $S$.

Again, if $T\simeq\pr{2}$, then $S\simeq\mathbb{F}_1$ and 
(P) holds for $S$;
hence we can suppose $T\not\simeq\pr{2}$. Let $E\subset S$ be the
exceptional curve, $p=\pi(E)\in T$,
$C\subset S$ an invariant curve different from $E$
and $B=\pi(C)\subset T$.
We can also assume $C^2\geq 0$, because if $C^2<0$ then $C$ is
extremal and there is nothing to prove.

Suppose $C\cdot E=0$, thus $C\cap E=\emptyset$. 
Since (P) is true for $T$, we have
\[ B\equiv\sum_i m_i C_i, \]
with $m_i\in\Z_{>0}$ and $C_i$ extremal in $\NE(T)$ for all $i$. 
Let $\widetilde{C_i}$ be the strict transform of $C_i$ in $S$:
$\widetilde{C_i}$ is extremal in $\NE(S)$ and
\[ \pi^*(C_i)=\widetilde{C_i}+\varepsilon_i E
\]
with $\varepsilon_i=0$ or $1$ if respectively $p\not\in C_i$ or $p\in
C_i$. Thus
\[ C=\pi^*(B)\equiv\sum_i m_i \pi^*(C_i)=\sum_i m_i \widetilde{C_i}+
(\sum_i m_i \varepsilon_i) E \]
is the decomposition we were looking for.

Now suppose $C\cdot E=1$. Let $C'$ be the other invariant curve in $S$
such that $C'\cdot E=1$, and let $B'=\pi(C')\subset T$.

In $\f_T$ we have: $B=V(x)$, $B'=V(x')$, $p=V(\langle x,x' \rangle)$. 
Let $y, y'\in G(\f_T)$ such that $\langle y, x \rangle\in\f_T$
and $\langle y', x' \rangle\in\f_T$. Then the
classes of $B$ and $B'$ in $\N(T)$ are respectively given by the relations:
\[ y+x'+ax=0,\qquad y'+x+a'x'=0  \]
with $a=B^2=C^2+1\geq 1$ and $a'=(B')^2$. 
Since $T\not\simeq\pr{2}$, there exists in
$\f_T$ a primitive collection of type $\{u,-u\}$: 
in fact, $T$ is obtained from some Hirzebruch surface $\mathbb{F}_a$ by
a finite sequence of smooth equivariant blow-ups, and all
Hirzebruch surfaces $\mathbb{F}_a$ contain a primitive collection of
type $\{u,-u\}$. 
The fact that $a\geq 1$ implies that the two cones $\langle x,x' \rangle$ and  $\langle y, x
\rangle$ cover strictly more than a half plane in $N_{\Q}$: therefore the only
possibility is $-x\in G(\f_T)$.
Since  $\langle y', x' \rangle\in\f_T$, $y'$ must lie in the cone 
$\langle -x,x'\rangle$, so $a'\leq 0$.

\vspace{10pt}

\hspace{100pt}\psfig{figure=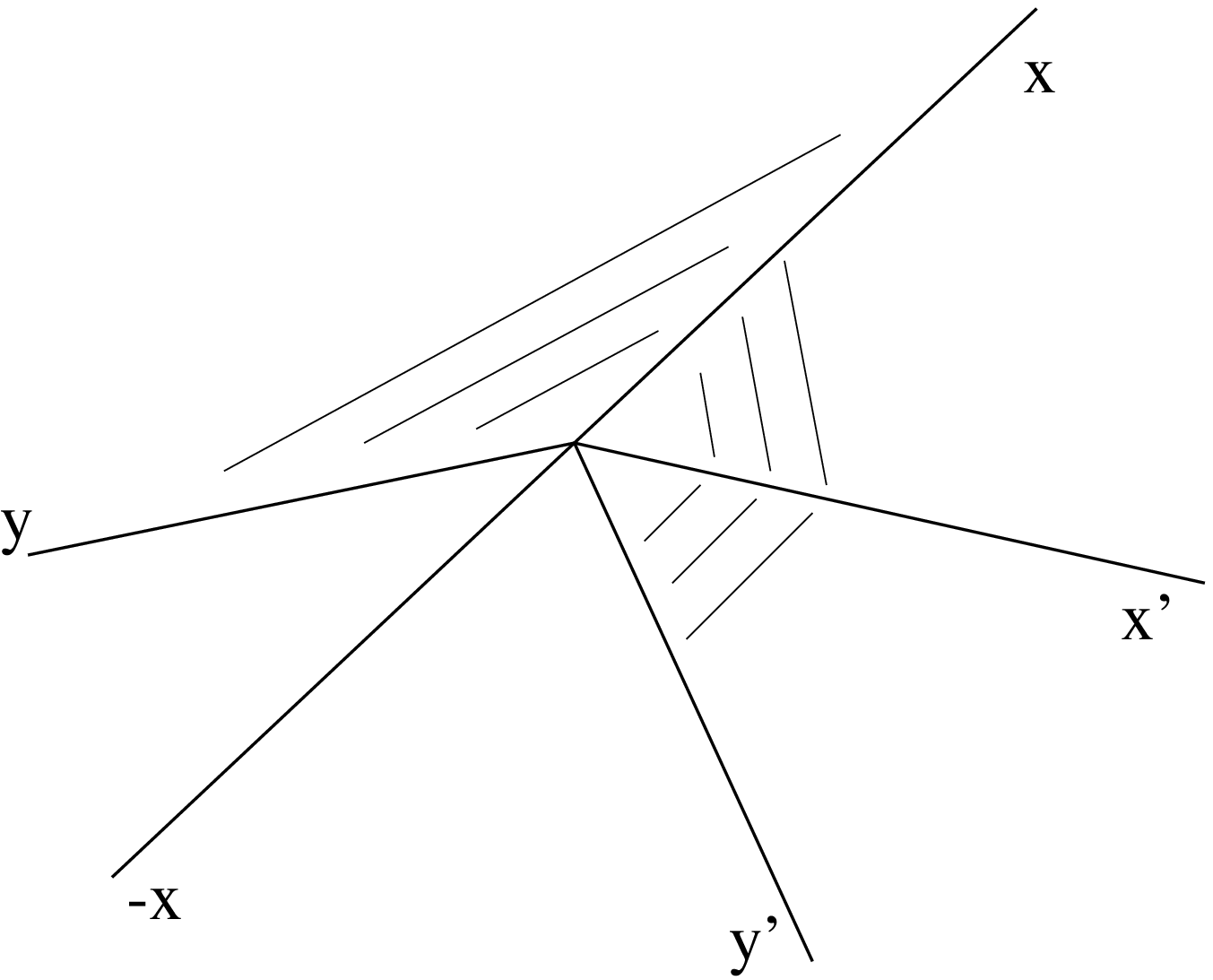,width=4cm} 
\begin{center}
{\footnotesize  $\f_{T}$}
\end{center}
\vspace{10pt}

Now we consider the element of $\A(T)$ given by the relation:
\[ y+(1-a')x'+(a-1)x - y'=0. \]
This class corresponds to $[B]-[B']$ and by lemma \ref{lemma} it is an
element of $\A(T)\cap\NE(T)$, because $1-a'>0$ and $a-1\geq 0$.
Thus we can apply (P) to  $[B]-[B']$ and get
\[ B\equiv B'+\sum_i m_i C_i, \]
$m_i\in\Z_{>0}$, $C_i$ extremal in $\NE(T)$ for all $i$.

In $S$ we have:
\begin{gather*}
 \pi^*(B)=C+E,\quad\pi^*(B')=C'+E,\quad\pi^*(C_i)=\widetilde{C_i}+
\varepsilon_i E
\quad\varepsilon_i=0,1 \\
\text{hence}\quad  
C\equiv C'+\sum_i m_i\widetilde{C_i}+ (\sum_i m_i\varepsilon_i)
E.
\end{gather*}
This is the decomposition we were looking for: $C'$ is
extremal, because $(C')^2=(B')^2-1=a'-1<0$.
\end{dimo}
\begin{dimo}[Proof of theorem \ref{bigteo}]
Let $H\in\Pic X$ be an ample line bundle and $C\subset X$ an 
invariant curve. 

If $C$ is contractible, there is nothing to prove. 
If $C$ is not
contractible, 
by theorem~\ref{guerra},
there exists an invariant surface $S$
containing an invariant curve $C'$ numerically equivalent to $C$, such
that $C'$ is not extremal in $\NE(S)$. 
In $S$ we have
\[ C'\equiv\sum_{i=1}^r m_i C_i, \]
with $r>1$, $m_i\in\Z_{>0}$, $C_i\subset S$. Thus in $X$
\[ C\equiv C'\equiv\sum_{i=1}^r m_i C_i. \]
Now $H\cdot C=\sum_{i=1}^r m_i H\cdot C_i$, so for
each $i$ we have $H\cdot C_i<H\cdot C$. If some of 
the $C_i$
are not contractible in $X$, we  iterate the procedure; since 
$H\cdot C$ is a positive integer, we must end after a finite
number of steps.
\end{dimo}

\emph{Question:} in the same assumptions of theorem~\ref{bigteo},
does every class in $\A(X)\cap\NE(X)$ decompose as a
linear combination with positive
integral coefficients of extremal classes?

Answering to this question means  understanding
 how a contractible non-extremal class $\eta$ decomposes as a linear
 combination of extremal classes. The point is that a priori
the coefficients are rational and can fail to be integers, 
so we don't know whether there exists an effective 
class $\gamma$ such that $\eta-\gamma$ is in $\NE(X)$. 

\medskip

\emph{Remark.} 
\label{rem}
In theorem \ref{bigteo}
the hypothesis of $X$ being projective is fundamental:
in fact, if $X$ is not projective, in general contractible classes 
can fail to generate the cone $\NE(X)$ even over $\Q$. 
As an example, consider the
non-projective toric variety $Y$ of dimension 3 and Picard number 4,
described on page~\pageref{esempio}.  
We recall that there is a smooth equivariant blow-up $\pi\colon X\rightarrow Y$
such that $X$ is projective.
In the fan of $Y$ there are seven primitive relations, four of which are contractible:
\[ \begin{array}{crclc}
\gamma_1:&\qquad e_1+f_2&=&f_1+e_2&\quad (C)\\
\omega_2:&\qquad f_1+e_3&=&e_1+f_3&\quad (C)\\
\omega_3:&\qquad e_2+f_3&=&f_2+e_3&\quad (C)\\
\gamma_4:&\qquad e_0+e_2&=&f_2& \\
\omega_5:&\qquad f_1+f_2+f_3&=&2e_0&\quad (C)\\
\gamma_6:&\qquad e_0+e_1&=&f_1& \\
\gamma_7:&\qquad e_0+e_3&=&f_3.&
\end{array} \]
We are keeping the notations of page~\pageref{esempio}.
The relations between these classes are:
\[ \gamma_1+\omega_2+\omega_3 =0,\quad \gamma_4=\omega_3+\gamma_7, 
\quad \gamma_6=\gamma_1+\gamma_4,\quad \gamma_7=\omega_2+\gamma_6. \]
These primitive relations generate all $\NE(Y)$, 
it is possible to see this looking at
$\pi_*\colon\NE(X)\rightarrow\NE(Y)$ (see corollary~\ref{basta}). 
Anyway, the three non-contractible primitive relations $\gamma_4$,
$\gamma_6$ and $\gamma_7$ can not be obtained as a linear combination
with positive rational coefficients of the four contractible
relations.  

We end this section with an application of theorem \ref{bigteo}:
\begin{prop}
\label{tie}
Suppose that $X$ is projective and let $H\subset X$ be an
ample divisor. Let 
\[ m_H=\min\{\, (C\cdot H)\,\,|\, C\text{ curve in }X\}. \] 
Then every class $\gamma\in\NE(X)$ such that $H\cdot
\gamma=m_H$ is extremal. 
\end{prop}
\begin{dimo}
By theorem \ref{bigteo}, we know that $\gamma$ is contractible; in
particular the relation associated to $\gamma$ is primitive:
\[ x_1+\cdots+x_r=a_1y_1+\cdots+a_ky_k. \]
Let $V_i=V(x_i)$ for $i=1,\ldo,r$ and consider the divisor
\[ \widetilde{H}=H-\frac{m_H}{r}(V_1+\cdots+V_r). \]
Since $V_i\cdot \gamma=1$ for all $i=1,\ldo,r$, we have
$\widetilde{H}\cdot \gamma=0$. 
Now let $B\subset X$ be an extremal curve, having numerical class
different from $\gamma$: 
let's show  that $\widetilde{H}\cdot B>0$.
Since $B$ is extremal, $B\cdot V_i\leq 1$ for all $i$:
we claim that it can not be $B\cdot V_i= 1$ for all $i$.
Indeed, in that case $[B]=r(P)$ with $P\supseteq\{x_1,\ldo,x_r\}$,
so $P=\{x_1,\ldo,x_r\}$ and $[B]=\gamma$.
Hence there exists an index $j\in\{1,\ldo,r\}$ such that $B\cdot
V_j\leq0$: this implies
\[ \widetilde{H}\cdot B\geq
H\cdot B-\frac{m_H}{r}\sum_{i\neq j}V_i\cdot B\geq 
m_H\left(1-\frac{r-1}{r}\right)>0. \]
\end{dimo}
\begin{cor}
Let $X$ be projective, $H\subset X$ an ample divisor and
$\gamma\in\NE(X)$ such that $\gamma\cdot H=1$. Then $\gamma$ is an
extremal class.
\end{cor}

\stepcounter{subsubsec}
\subsubsection{Blow-ups}
Throughout this section, $X$ and $Y$ are smooth, complete,
$n$-dimensional toric varieties and $\pi\colon X\rightarrow
Y$ is a smooth equivariant blow-up.  
We want to study the behaviour of contractible
classes (and more generally of primitive relations)
under $\pi$. 
We recall that any birational map between two smooth, complete 
toric varieties factorizes as a sequence of smooth equivariant blow-ups 
and blow-downs \cite{morelli,abr1}, hence it is natural to study 
what happens through a single step.

Let $A=V(\langle 
v_1,\ldo,v_r\rangle)\subset Y$ be the center of the blow-up,
 $E=V(v)\subset X$  the
exceptional divisor and  $\delta\in\NE(X)$ the contractible class
coming from the blow-up, with relation $v_1+\cdots +v_r=v$. 

We recall that $G(\fx)=G(\fy)\cup\{v\}$ and 
$\N(Y)=\N(X)/\Q\cdot\delta$; actually $\N(Y)$ is naturally 
identified with the hyperplane in $\N(X)$ given by those relations
where $v$ does not appear. 
 Under the projection
$\pi_*\colon\N(X)\rightarrow\N(Y)$, we have $\NE(Y)=\pi_*(\NE(X))$.

We recall that when $X$ and $Y$ are projective, the class $\delta$ is
extremal in $\NE(X)$ and $\NE(Y)$ is obtained projecting $\NE(X)$ from 
$\delta$. Hence every extremal ray of $\NE(Y)$ comes from an extremal
ray of $\NE(X)$. Anyway, since $\NE(X)$ in general is not a simplicial
cone, it may happen that some extremal ray of $\NE(X)$ is projected on
something non-extremal in $\NE(Y)$. It is natural to ask if similar
properties hold for contractible classes. The result we get is the
following: 
\begin{teo}
\label{piove}
Let $\mathcal{C}_X$ (respectively, $\mathcal{C}_Y$) be the set of
contractible classes in $\NE(X)$ (respectively, $\NE(Y)$). Then:
\begin{itemize}
\item 
$\pi_*\colon\NE(X)\rightarrow\NE(Y)$ is iniective on $\mathcal{C}_X\smallsetminus\{\delta\}$;
\item $\pi_*(\mathcal{C}_X\smallsetminus\{\delta\})$ are all primitive
  relations;
\item
  $\mathcal{C}_Y\subseteq\pi_*(\mathcal{C}_X\smallsetminus\{\delta\})$. 
\end{itemize}
Moreover, if $\omega\in\mathcal{C}_X\smallsetminus\{\delta\}$ 
and $\pi_*(\omega)$ is not contractible, then $\omega\cdot E=1$ and
there exists some $i\in\{1,\ldo,r\}$ such that $\omega\cdot V(v_i)=-1$.
\end{teo}
\noindent\emph{Remarks:}

\noindent $\bullet\ \,$ in general $\mathcal{C}_Y\subsetneqq
\pi_*(\mathcal{C}_X\smallsetminus\{\delta\})$, namely 
 $\omega\in\NE(X)$ contractible does not 
imply  $\pi_*(\omega)$ contractible in $\NE(Y)$. 
We can see this in dimension 2 for extremal
classes: 
let $S$ be a smooth toric surface containing an invariant curve $C$,
not extremal,  such
that $C^2=0$ (for instance, take the blow-up of $\pr{1}\times\pr{1}$
at one point). Consider the blow-up of $S$ at
a fixed point of $C$. The strict transform $\widetilde{C}$ has
self-intersection $-1$, hence it is extremal, but its class projects
on the class of $C$, which is not extremal. 

\noindent $\bullet\ \,$  Suppose that $X$, $Y$ are projective and
let $\gamma\in\NE(Y)$ and $\widetilde{\gamma}\in\NE(X)$
be two contractible classes such that
$\pi_*(\widetilde{\gamma})=\gamma$.
Then we have:
\[ \gamma \text{ extremal}\quad\Longrightarrow\quad
\widetilde{\gamma}\text{ extremal, }\]
because there is always an extremal class that projects onto $\gamma$,
but since there is a unique contractible class projecting onto $\gamma$,
this extremal class must be $\widetilde{\gamma}$.
The converse is false in general: even if $\widetilde{\gamma}$ is
extremal, in $\NE(X)$ there can be a relation 
$a\widetilde{\gamma}+b\delta=\sum_i\lambda_i\omega_i$, with
$a,b,\lambda_i\in\mathbb{Z}_{>0}$ and $\omega_i$ extremals. Then
$a\gamma=\sum_i\lambda_i\pi_*(\omega_i)$ is not extremal. We remark
that in this way, if $a>1$, we get rational coefficients.

\smallskip

There is a result of H.~Sato that allows 
to compute the primitive collections of
$X$ from the ones of $Y$:
\begin{teo}[H.~Sato \cite{sato}, Theorem 4.3]
\label{sato2}
Let $P\in\PC(\fy)$.
\begin{enumerate}[a.]
\item If
$P\cap\{v_1,\ldo,v_r\}=\emptyset$, then $P$ is a primitive
collection in $\fx$.
\item If
$P\supset\{v_1,\ldo,v_r\}$, then
$(P\smallsetminus\{v_1,\ldo,v_r\})\cup\{v\}$ is a primitive
collection in $\fx$.
\item If
$P\cap\{v_1,\ldo,v_r\}\neq\emptyset$ and
$P\not\supset\{v_1,\ldots,v_r\}$, then $P$ is a primitive collection
in $\fx$; moreover, either 
$(P\smallsetminus\{v_1,\ldo,v_r\})\cup\{v\}$ is a primitive
collection in $\fx$, or it contains a primitive collection
$Q'\in\PC(\fx)$ such that
$Q'=(Q\smallsetminus\{v_1,\ldo,v_r\})\cup\{v\}$ for some $Q\in\PC(\fy)$.
\end{enumerate}

Except $\{v_1,\ldo,v_r\}$, all the primitive collections of $\fx$ are
obtained in one of these three ways.
\end{teo}
\begin{defi} We  say that a primitive collection $P$
is of type $a$, $b$ or $c$ if respectively 
$P\cap\{v_1,\ldo,v_r\}=\emptyset$, $P\supset\{v_1,\ldo,v_r\}$ or 
$P\cap\{v_1,\ldo,v_r\}=\{v_{i_1},\ldo,v_{i_m}\}$ with $0<m<r$.
Moreover, for each $P$ we define $P'$ in the following way:
$P'=P$ is $P$ is of type $a$; $P'=(P\smallsetminus\{v_1,\ldo,v_r\})\cup\{v\}$ 
if $P$ is of type $b$ or $c$. According to theorem~\ref{sato2},
$P'$ is always a primitive collection in $\fx$ when $P$ is of type $a$
or $b$. 
\end{defi}

We are going to 
analyze, for each of these types of primitive collections, 
how primitive relations change under the blow-up.
Since the same $P$ can be a primitive collection for both $\fx$ and
$\fy$, we will denote by $r_X(P)$ and
$r_Y(P)$ respectively the primitive relations associated to $P$ in
$\fx$ and $\fy$.

We remark that even if these three types are clearly distinct from a
combinatorial point of view, this distinction doesn't have a clear
geometrical meaning. The relative positions of $A$ and the locus $Z$ of
$r_Y(P)$ can be quite different inside each type $a$, $b$, $c$.
For instance, when $P$ is of type $a$ and $r_Y(P)$ is contractible, it
can be either $Z\cap A=\emptyset$ or $Z\cap A\neq\emptyset$, and even
$Z\supseteq A$ or $Z\subseteq A$.

\begin{lemma}
\label{a}
Let $P\in\PC(\fy)$ be of type $a$ or $b$.
Then $\pi_*(r_X(P'))=r_Y(P)$, and $r_X(P')$ is contractible
if and only if $r_Y(P)$ is contractible.
\end{lemma}
For primitive relations of type $c$, 
the analysis is more delicate:
theorem~\ref{sato2} tells us that $P$ is a primitive
collection for $\fx$, and
$P'$ may also be a
primitive collection for $\fx$. It is easy to see that that 
the primitive
relation of $P$ in $\fx$ is the same that the one in $\fy$:
$r_X(P)=r_Y(P)$.
 But if $P'$ is primitive, 
the primitive relation $r_X(P')$ is unknown, and
in general $\pi_*(r_X(P'))$ is not a primitive relation. 
\begin{lemma}
\label{c}
Let $P\in\PC(\fy)$ be of type $c$, 
with associated  primitive relation $r_Y(P)$:
\[
v_1+\cdots+v_m+x_1+\cdots+x_h=a_{m+1}v_{m+1}+\cdots+a_tv_t+b_1y_1+\cdots
+b_ky_k
\]
with $1\leq m<r$, $m\leq t\leq r$, and let
$Z=V(\langle v_{m+1},\ld,v_t,y_1,\ld,y_k\rangle)$ be the
locus of $r_Y(P)$. 
Then:
\begin{enumerate}[(i)]
\item
$r_X(P)$ is contractible if and only if $r_Y(P)$ is contractible
and $Z\cap A=\emptyset$;
in this case, even if $P'$ is primitive, 
$\pi_*(r_X(P'))\neq r_Y(P)$;
\item
if $r_Y(P)$ is contractible and $Z\cap A\neq\emptyset$, then $P'$ is a
primitive collection in $\fx$,
$\pi_*(r_X(P'))=r_Y(P)$  and $r_X(P')$ is contractible;
\item
if $P'$ is a primitive collection in $\fx$ and
$r_X(P')$ is contractible, then $\pi_*(r_X(P'))=r_Y(P)$; if
moreover  $r_X(P')\cdot V(v_i)\neq
-1$ for all $i=1,\ldo,r$, then $r_Y(P)$ is contractible.
\end{enumerate}
\end{lemma}
\noindent The proofs of lemma \ref{a} and lemma \ref{c}
are completely combinatorial and we postpone them to the 
end of the section.

\begin{dimo}[Proof of theorem \ref{piove}]
Let $\omega=r_X(Q)\in\NE(X)$ be a contractible class,
$\omega\neq\delta$.
Let's show that $\pi_*(\omega)$ is a primitive relation. 
By theorem~\ref{sato2}, either $Q=P'$ for some primitive collection $P$
in $\fy$ of type $a$, $b$ or $c$, or $Q$ is a primitive
collection also in $\fy$ and has type $c$.

If $Q=P'$ with $P$ of type $a$ or $b$, then by lemma~\ref{a}
$\pi_*(\omega)=r_Y(P)$ and $r_Y(P)$ is contractible.

If $Q$ is a primitive collection also in $\fy$ and has type $c$, then
by lemma~\ref{c} $(i)$, $\pi_*(\omega)=r_Y(Q)$ and $r_Y(Q)$ is
contractible.

If $Q=P'$ with $P$ of type $c$,  then by lemma~\ref{c}~$(iii)$,
$\pi_*(\omega)=r_Y(P)$ is primitive. Moreover, if $\pi_*(\omega)$ is
not contractible, it must be $\omega\cdot V(v_i)=-1$ for some
$i\in\{1, \ldo,r\}$. Since $v\in Q$, we have $\omega\cdot E=1$.

The only possible case where a primitive collection $P$ of $\fy$ gives
two primitive collections $P$, $P'$ in $\fx$ is when $P$ is of type
$c$; in this case, by lemma~\ref{c},
the two relations $r_X(P)$ and $r_X(P')$ can not be
both contractible.
Hence $\pi_*$ is injective on $\mathcal{C}_X\smallsetminus\{\delta\}$.

Let's show that
$\mathcal{C}_Y\subseteq\pi_*(\mathcal{C}_X\smallsetminus\{\delta\})$. 
Let $\gamma\in\NE(Y)$ be a contractible class; then $\gamma=r_Y(P)$
for a primitive collection $P$ in $\fy$. 

If $P$ is of type $a$ or $b$, then $P'$ is a primitive collection in
$\fx$ by theorem~\ref{sato2}. By lemma~\ref{a}, $r_X(P')$ is
contractible and  $\pi_*(r_X(P'))=r_Y(P)$. 

Suppose that $P$ is of type $c$; then by theorem~\ref{sato2} $P$ is a
 primitive collection in $\fx$ too. 
 Let $Z\subset Y$ be the locus of
$\gamma$. If $Z\cap A=\emptyset$, then lemma~\ref{c}~$(i)$ applies,
so  $r_X(P)$ is contractible and $\pi_*(r_X(P))=r_Y(P)$.

If $Z\cap A\neq \emptyset$, lemma~\ref{c}~$(ii)$ applies, so 
$P'$ is a primitive collection, $r_X(P')$ is contractible 
and $\pi_*(r_X(P'))=r_Y(P)$. 
\end{dimo}

\begin{cor}
\label{basta}
Let $Y$ be a complete, non-projective, smooth toric variety such that there
exists a smooth equivariant blow-up $X\rightarrow Y$ with $X$
projective. Then every class $\eta\in\A(Y)\cap\NE(Y)$ is a linear
combination with positive integral coefficients of primitive
classes. In particular,
\[ \NE(Y)=\sum_{P\in\PC(\fy)} \Q_{\,\geq 0}r(P). \]
\end{cor}
\begin{dimo}
It is an immediate consequence of theorems \ref{bigteo} and \ref{piove}.
\end{dimo}
\begin{cor}
Let $\gamma=r_Y(P)$ be a contractible class in $\NE(Y)$ and 
$\widetilde{\gamma}\in\NE(X)$ contractible
such that $\pi_*(\widetilde{\gamma})=\gamma$. 
Then we get a commutative diagram:
\[
\xymatrix{
 X \ar[r]^{\pi} \ar[d]^{\ph_{\widetilde{\gamma}}} &  
{Y} \ar[d]^{\ph_{\gamma}}  \\
{X_{\widetilde{\gamma}}} \ar[r]_{\psi} & {Y_{\gamma}} }
\]
If $\ph_{\gamma}$ is a blow-up and $P$ is of type $a$ or $b$, 
then all the four morphisms are blow-ups. 

Suppose that $\ph_{\gamma}$ is a $\pr{h-1}$-bundle.
If $P$ is of type $a$, then 
also $\ph_{\widetilde{\gamma}}$ is a $\pr{h-1}$-bundle; 
$A$ is a $\pr{h-1}$-bundle on
$\ph_{\gamma}(A)$ and $\psi$ is the blow-up of
$Y_{\gamma}$  along $\ph_{\gamma}(A)$.
 If $P$ is of type $b$, then 
$\ph_{\widetilde{\gamma}}$ is a $\pr{h-r}$-bundle, and $\psi$ is 
a $\pr{r-1}$-bundle.
\end{cor}

\noindent \emph{Some examples with $\gamma$ numerically effective:}

\noindent $a)$ let $Y=\mathbb{P}_{\pr{2}}(\mathcal{O}\oplus\mathcal{O}(1))$,
$\ph_{\gamma}$ the fibration on
$\pr{2}$, and let $\pi\colon X\rightarrow Y$ be the blow-up along the
fiber of a fixed point of  $\pr{2}$.
 We get:
\[
\xymatrix{
X \ar[r]
^{\pi\qquad\ \ } 
\ar[d]^{\ph_{\widetilde{\gamma}}} &  
{\mathbb{P}_{\pr{2}}(\mathcal{O}\oplus\mathcal{O}(1))} 
\ar[d]^{\ph_{\gamma}}  \\
{\mathbb{F}_1} \ar[r]_{\psi} & {\pr{2}} }
\]

\noindent $b)$ 
let $X=\mathbb{F}_1\times \pr{1}$, $Y=\pr{2}\times \pr{1}$, $\pi$ the
blow-up along $\{*\}\times\pr{1}$ and $\ph_{\gamma}$ the projection on
$\pr{1}$. We get:
\[
\xymatrix{
{\mathbb{F}_1\times\pr{1}} \ar[r]^{\pi} \ar[d]^{\ph_{\widetilde{\gamma}}} &  
{\pr{2}\times\pr{1}} \ar[d]^{\ph_{\gamma}}  \\
{\pr{1}\times\pr{1}} \ar[r]_{\psi} & {\pr{1}} }
\]

\noindent $c)$ 
let $Y=\pr{2}\times \pr{2}$, $\pi\colon X\rightarrow Y$ the
blow-up at one fixed point $p$ and $\ph_{\gamma}$ one of the projections on
$\pr{2}$. Let $F$ be the $\pr{2}$ over $\ph_{\gamma}(p)$ in $Y$:  
$F$ is blown-up to a surface 
$\mathbb{F}_1$ in $X$, that intersect the exceptional divisor along 
its exceptional curve $L$; $\widetilde{\gamma}$ is the 
numerical class 
of the proper tranform of a line through $p$ in $F$. 
The morphism $\ph_{\widetilde{\gamma}}$ is birational: the exceptional 
locus is 
the surface $\mathbb{F}_1$, that is contracted on $L$. 
The image $X_{\widetilde{\gamma}}$ is singular and $\psi$ has as 
a general fiber a $\pr{2}$, while the fiber over $p$ has dimension 3.
\[
\xymatrix{
{X} \ar[r]^{\pi\ \ \ } \ar[d]^{\ph_{\widetilde{\gamma}}} &  
{\pr{2}\times\pr{2}} \ar[d]^{\ph_{\gamma}}  \\
{X_{\widetilde{\gamma}}} \ar[r]_{\psi} & {\pr{2}} }
\]

\begin{dimo}[Proof of lemma \ref{a}]

\noindent \emph{Case a: $P\cap\{v_1,\ld,v_r\}=\emptyset$.}
Consider the primitive relation $r_Y(P)$:
\[ x_1+\cdots +x_h=a_1y_1+\cdots+a_ky_k. \]
We suppose first that $\langle y_1,\ld,y_k\rangle\not\supseteq\langle
v_1,\ld,v_r\rangle$: then $\langle y_1,\ld,y_k\rangle\in\fx$ and
$r_X(P)=r_Y(P)$. 

Suppose that $r_X(P)$ is contractible and let $\nu=\langle z_1,\ld,z_s
\rangle\in\fy$ such that
$\{z_1,\ld,z_s\}\cap\{x_1,\ld,x_h,y_1,\ld,y_k\}=\emptyset$ and
$\langle y_1,\ldots,y_k\rangle+\nu\in\fy$. If $\langle
y_1,\ld,y_k\rangle+\nu$ does not contain  
$\langle v_1,\ld,v_r\rangle$, then it is a cone in $\fx$ too, and
it is clear by the contractibility of 
$r_X(P)$ that $\langle
x_1,\ld,\check{x}_i,\ld,x_h,y_1,\ld,y_k\rangle+\nu\in\fy$ for all
$i=1,\ld,h$. If 
$\langle y_1,\ld,y_k\rangle+\nu$ contains $\langle v_1,\ld,v_r\rangle$, 
then the set $\{y_1,\ld,y_k,z_1,\ld,z_s,v\}\smallsetminus\{v_j\}$ generates
a cone in $\fx$ for all $j=1,\ld,r$. We can choose $j$ such that
$v_j\not\in\{y_1,\ld,y_k\}$: then 
for the contractibility of $r_X(P)$ 
the set
$\{x_1,\ld,x_h,y_1,\ld,y_k,z_1,\ld,z_s,v\}\smallsetminus\{x_i,v_j\}$
generates a cone in $\fx$ for all $i=1,\ld,h$.
Hence 
$\langle x_1,\ld,\check{x}_i,\ld,x_h,
y_1,\ld,y_k,z_1,\ld,z_s\rangle\in\fy$ for all $i=1,\ld,h$.
So $r_Y(P)$ is
contractible.

Conversely, suppose that $r_Y(P)$ is contractible and let $\eta
=\langle z_1,\ld,z_s\rangle\in\fx$
such that $\{z_1,\ld,z_s\}\cap\{ x_1,\ld,x_h,y_1,\ld,y_k\}= 
\emptyset$ and
$\langle y_1,\ld,y_k\rangle+\eta\in\fx$. If $v\not\in\eta$, it is
easy to see that $\langle
x_1,\ld,\check{x}_i,\ld,x_h,y_1,\ld,y_k\rangle+\eta\in\fx$ for all
$i=1,\ld,h$. If $v=z_1$, then
$\langle y_1,\ld,y_k,
v_1,\ld,v_r,z_2,\ld,z_s \rangle\in\fy$, so by the
contractibility of $r_Y(P)$  $\langle
x_1,\ld,\check{x}_i,\ld,x_h, y_1,\ld,y_k,v,z_2,\ld,z_s\rangle\in\fx$
for all $i=1,\ld,h$. Hence $r_X(P)$ is contractible.

We have now to consider the case where the primitive relation $r_Y(P)$
is
\[ x_1+\cdots+x_h=a_1y_1+\cdots+a_ky_k+b_1v_1+\cdots +b_rv_r. \]
We order the $v_i$ in such a way that
$b_1=\cdots=b_{p-1}<b_p\leq b_{p+1}\leq \cdots\leq
b_r$, with $p=2,\ld,r+1$. 
Consider the relation $r_Y(P)+b_1\delta$:
\[ x_1+\cdots+x_h=a_1y_1+\cdots+a_ky_k+b_1v+(b_p-b_1)v_p+\cdots +(b_r-b_1)v_r. \]
Since $\langle y_1,\ld,y_k,v_1,\ld,v_r\rangle\in\fy$, we have
$\langle y_1,\ld,y_k,v,v_p,\ld,v_r\rangle\in\fx$: so 
$r_X(P)=r_Y(P)+b_1\delta$ and $\pi_*(r_X(P))=r_Y(P)$. 

Suppose that $r_X(P)$ is contractible and let $\nu=\langle z_1,\ld,z_s
\rangle\in\fy$ such that
$\{z_1,\ld,z_s\}\cap\{ 
x_1,\ld,x_h,y_1,\ld,y_k,v_1,\ld,v_r\}=\emptyset$ 
and $\langle y_1,\ld,y_k,v_1,\ld,v_r\rangle+\nu\in\fy$. Then
$\langle y_1,\ld,y_k,v,v_p,\ld,v_r\rangle+\nu$ is in $\fx$ and by
the contractibility of $r_X(P)$, we get $\langle x_1,\ld,\check{x}_i,\ld,x_h,
y_1,\ld,y_k,v_1,\ld,v_r \rangle+\nu\in\fy$ 
for all $i=1,\ld,h$. So $r_Y(P)$ is
contractible.

Conversely, suppose that $r_Y(P)$ is contractible and let 
$\eta=\langle z_1,\ld,z_s\rangle\in\fx$
such that $\{z_1,\ld,z_s\}\cap\{ 
x_1,\ld,x_h,y_1,\ld,y_k,v,v_p,\ld,v_r\}=\emptyset$ and
\[ \langle y_1,\ld,y_k,v,v_p,\ld,v_r\rangle+\eta\in\fx.\] 
Then $\langle y_1,\ld,y_k,v_1,\ld,v_r\rangle+\eta$ is in $\fy$ and 
the contractibility of $r_Y(P)$ implies  
$\langle x_1,\ld,\check{x}_i,\ld,x_h,
y_1,\ld,y_k,v,v_p,\ld,v_r \rangle+\eta\in\fx$ for all $i=1,\ld,h$. 
Hence $r_X(P)$ is
contractible.

\medskip

\noindent \emph{Case b: $P\supset\{v_1,\ld,v_r\}$}.
Let the primitive relation $r_Y(P)$ be \[
v_1+\cdots+v_r+x_1+\cdots+x_h=a_1y_1+\cdots+a_ky_k. \]
Clearly $\langle y_1,\ld,y_k\rangle$ is also in $\fx$ and the
primitive relation associated to $P'$ in $\fx$ is
$v+x_1+\cdots+x_h=a_1y_1+\cdots+a_ky_k$.
Thus $r_X(P')=r_Y(P)-\delta$ and $\pi_*(r_X(P'))=r_Y(P)$. 

It is easy to see that $r_X(P')$ is contractible if and only if
$r_Y(P)$ is contractible.
\end{dimo}
\begin{dimo}[Proof of lemma \ref{c}]
We recall that the primitive relation $r_Y(P)$ is
\[
v_1+\cdots+v_m+x_1+\cdots+x_h=a_{m+1}v_{m+1}+\cdots+a_tv_t+b_1y_1+\cdots
+b_ky_k
\]
with $1\leq m<r$, $m\leq t\leq r$. 

$(i)$ Suppose that $r_X(P)$ is contractible, and
let $\nu=\langle z_1,\ld,z_s\rangle\in\fy$ such that
$\{z_1,\ld,z_s\}\cap\{
v_1,\ld,v_t,x_1,\ld,x_h,y_1,\ld,y_k\}=\emptyset$ and 
\[\langle v_{m+1},\ld,v_t,y_1,\ld,y_k\rangle+\nu\in\fy.\]
Then the cone 
$\langle v_{m+1},\ld,v_t,y_1,\ld,y_k\rangle+\nu$ is also in $\fx$, and
the contractibility of $r_X(P)$ implies that $\langle
v_1,\ld,v_t,x_1,\ld,\check{x}_i,\ld,x_h,y_1,\ld,y_k\rangle+\nu\in\fy$ for
all $i=1,\ld,h$ and $\langle
v_1,\ld,\check{v}_j,\ld,v_t,x_1,\ld,x_h,y_1,\ld,y_k\rangle+\nu\in\fy$ for
all $j=1,\ld,m$. So $r_Y(P)$ is contractible.

To show that $Z\cap A=\emptyset$, we have to show that the cone
$\langle y_1,\ld,y_k,v_1,\ld,v_r\rangle$ is not in $\fy$. If this cone
were in $\fy$, we would get $\langle y_1,\ld,y_k,v\rangle\in\fx$,
hence $\langle x_1,\ld,x_h,v\rangle\in\fx$ which implies $\langle
P\rangle \in \fx$, a contradiction.

Conversely, let's suppose that $r_Y(P)$ is contractible and that 
$Z\cap A=\emptyset$. Hence we have $\langle
y_1,\ld,y_k,v_1,\ld,v_r \rangle\not\in\fy$, that implies  $\langle
y_1,\ld,y_k,v \rangle\not\in\fx$. 
Let's show that $r_X(P)$ is contractible.
Let $\eta=\langle z_1,\ld,z_s\rangle$ be such that 
$\{z_1,\ld,z_s\}\cap\{
v_1,\ld,v_t,x_1,\ld,x_h,y_1,\ld,y_k\}=\emptyset$ and 
$\langle v_{m+1},\ld,v_t,y_1,\ld,y_k\rangle+\eta\in\fx$. Then
$v\not\in\eta$ and by the contractibility of $r_Y(P)$ we get that the
cones  $\langle
v_1,\ld,v_t,x_1,\ld,\check{x}_i,\ld,x_h,y_1,\ld,y_k\rangle+\eta$, 
$\langle
x_1,\ld,x_h,v_1,\ld,\check{v}_j,\ld,v_t,y_1,\ld,y_k\rangle+\eta$
are in $\fx$ for all $i=1,\ld,h$ and $j=1,\ld,m$. 
Thus $r_X(P)$ is contractible. 

Finally, suppose that $P'=\{v,x_1,\ld,x_h\}$ is primitive. 
We remark that  $r_Y(P)-\delta$ is 
\begin{multline*}
v+x_1+\cdots+x_h= \\
(a_{m+1}+1)v_{m+1}+\cdots+(a_t+1)v_t+v_{t+1}+\cdots+
v_r+b_1y_1+\cdots+b_ky_k.
\end{multline*}
In general $r_Y(P)-\delta$ fails to be the primitive relation
associated to $P'$ in $\fx$, because the cone $\langle
v_{m+1},\ld,v_r,y_1,\ld,y_k \rangle$ does not have to belong to $\fx$.
We have:
\begin{eqnarray*}
\pi_*(r_X(P'))= r_Y(P)\ \Longleftrightarrow\
r_X(P')=r_Y(P)+\lambda\delta,\ \lambda\in\Q\ \Longleftrightarrow\\
\ r_X(P')=r_Y(P)-\delta\ \Longleftrightarrow\ 
\langle
v_{m+1},\ld,v_r,y_1,\ld,y_k \rangle\in\fx. 
\end{eqnarray*} 
Under our hypotheses, 
since  $\langle
v_1,\ld,v_r,y_1,\ld,y_k \rangle\not\in\fy$ and $r_Y(P)$ is
contractible, the cone  $\langle
v_{m+1},\ld,v_r,y_1,\ld,y_k \rangle$ can not be in $\fy$, 
neither in $\fx$. Therefore $\pi_*(r_X(P'))\neq r_Y(P)$.

\medskip

$(ii)$ Assume that $r_Y(P)$ is contractible and that 
 $Z\cap A\neq\emptyset$, namely  that the cone
$\langle
v_1,\ld,v_r,y_1,\ld,y_k \rangle$ is in $\fy$.
Hence we get 
\[\langle
v_1,\ld,v_r,x_1,\ld,\check{x}_i,\ld,x_h\rangle\in\fy\text{
 for all }i=1,\ld,h;  \] moreover $\langle
x_1,\ld,x_h\rangle\in\fy$. 
Therefore $\langle
x_1,\ld,x_h\rangle$ and
$\langle
v,x_1,\ld,\check{x}_i,\ld,x_h\rangle$ are in $\fx$ for all
$i=1,\ld,h$. Since $\langle
v,x_1,\ld,x_h\rangle\not\in\fx$, $P'=\{v,x_1,\ld,x_h\}$ is a
primitive collection in $\fx$. Moreover $\langle
v_{m+1},\ld,v_r,y_1,\ld,y_k \rangle\in\fx$ implies that
$r_X(P')=r_Y(P)-\delta$ and $\pi_*(r_X(P'))=r_Y(P)$.

Let's show that $r_X(P')$ is contractible: 
we consider $\eta=\langle z_1,\ld,z_s\rangle\in\fx$ such that 
$\{z_1,\ld,z_s\}\cap
\{
v,v_{m+1},\ld,v_r,x_1,\ld,x_h,y_1,\ld,y_k\}=\emptyset$ and 
\[\langle
v_{m+1},\ld,v_r,y_1,\ld,y_k\rangle+\eta\in\fx.\] 
Since $v\not\in\eta$,  $\langle
v_{m+1},\ld,v_r,y_1,\ld,y_k\rangle+\eta\in\fy$ and the contractibility of
$r_Y(P)$ implies  $\langle
v_1,\ld,\check{v}_j,\ld,v_r,x_1,\ld,x_h,y_1,\ld,y_k\rangle+ 
\eta\in\fy$ for all $j=1,\ld,m$, and  $\langle
v_1,\ld,v_r,x_1,\ld,\check{x}_i,\ld,x_h,y_1,\ld,y_k\rangle+\eta\in\fy$ 
for all $i=1,\ld,h$. Hence in $X$ we get  
that the cones
$\langle
v_{m+1},\ld,v_r,x_1,\ld,x_h,y_1,\ld,y_k\rangle+\eta$ and  
$\langle v,v_{m+1},\ld,v_r,x_1,\ld,\check{x}_i,\ld,x_h,y_1,\ld,y_k 
\rangle+\eta$ are in $\fx$
for all $i=1,\ld,h$. So $r_X(P')$ is contractible.   

\medskip

$(iii)$ Assume that $P'=\{v,x_1,\ld,x_h\}$ is primitive,
and that $r_X(P')$ is contractible. We have to show that
$r_X(P')+\delta=r_Y(P)$. Let $r_X(P')$ be:
\[ v+x_1+\cdots+x_h=c_{i_1} v_{i_1}+\cdots +c_{i_p}v_{i_p}+d_1w_1+
\cdots+d_lw_l, \]
where $p<r$ and $l\geq 0$. 
We order the $c_i$ in such a way that $c_{i_1}\geq \cdots\geq 
c_{i_s}>1=c_{i_{s+1}}=\cdots= c_{i_p}$, with $s=0,\ld,p$, 
and set $\{v_{j_1},\ld,v_{j_q}\}= 
\{v_1,\ld,v_r\}\smallsetminus\{v_{i_1},\ld,v_{i_p}\}$.

Since $r_X(P')$ is contractible, we have
$\langle v,w_1,\ld,w_l\rangle\in\fx$, so 
\[\langle
v,v_1,\ld,\check{v}_j,\ld,v_r,w_1,\ld,w_l\rangle\in\fx\ \text{ for all }
j=1,\ld,r. \]
Hence for all
$e=1,\ld,q$ we get $\langle
x_1,\ld,x_h,v_1,\ld,\check{v}_{j_e},\ld,v_r\rangle\in\fx$.  
This implies that $P$ is not contained in the set
$\{x_1,\ld,x_h,v_1,\ld,\check{v}_{j_e},\ld,v_r\}$ for all
$e=1,\ld,q$, so $P\supseteq\{x_1,\ld,x_h,v_{j_1},\ld,v_{j_q}\}$. 

On the other hand, we have the relation $r_X(P')+\delta$:
\begin{multline*}
v_{j_1}+\cdots +v_{j_q}+x_1+\cdots+x_h=\\
(c_{i_1}-1) v_{i_1}+\cdots 
+(c_{i_s}-1)v_{i_s}+d_1w_1+\cdots+d_lw_l. \end{multline*}
Since the cone $\langle
v_{i_1},\ld,v_{i_s},w_1,\ld,w_l\rangle$
is in $\fx$, the set $\{x_1,\ld,x_h,v_{j_1},\ld,v_{j_q}\}$ can not
generate a cone in $\fx$: therefore
$P=\{x_1,\ld,x_h,v_{j_1},\ld,v_{j_q}\}$, and $P$ is the unique primitive
collection of $\fy$ contained in
$\{x_1,\ld,x_h,v_1,\ld,v_r\}$. Moreover,
$r_X(P')+\delta=r_Y(P)$, so $\pi_*(r_X(P'))=r_Y(P)$.

We suppose now that
$r_X(P')\cdot V(v_i)\neq
-1$ for all $i=1,\ld,r$. Since $r_X(P')=r_Y(P)-\delta$, it is
\begin{multline*} 
v+x_1+\cdots+x_h= \\
(a_{m+1}+1)v_{m+1}+\cdots+(a_t+1)v_t+v_{t+1}+\cdots+
v_r+b_1y_1+\cdots+b_ky_k.
\end{multline*}
Therefore $r_X(P')\cdot V(v_i)\neq
-1$ for all $i=1,\ld,r$ is equivalent to $t=r$, namely $r_Y(P)$ is
\[
v_1+\cdots+v_m+x_1+\cdots+x_h=a_{m+1}v_{m+1}+\cdots+a_rv_r+b_1y_1+\cdots
+b_ky_k.
\]
We want to show that $r_Y(P)$
is contractible. Let $\nu=\langle z_1,\ld,z_s\rangle
\in\fy$ such that
$\{z_1,\ld,z_s\}\cap\{v_1,\ld,v_r,x_1,\ld,x_h,y_1,\ld,y_k\}=\emptyset$
and
\[\langle v_{m+1},\ld,v_r,y_1,\ld,y_k\rangle+\nu\in\fy. \]
Then $\langle
v_{m+1},\ld,v_r,y_1,\ld,y_k\rangle+\nu\in\fx$ and the contractibility of
$r_X(P')$ and $\delta$ implies $\langle
v,v_1,\ld,\check{v}_j,\ld,v_r,x_1,\ld,\check{x}_i,
\ld,x_h,y_1,\ld,y_k\rangle+\nu\in\fx$   
for all $j=1,\ld,r$ and $i=1,\ld,h$. 
Hence also $\langle
v_1,\ld,\check{v}_j,\ld,v_r,x_1,\ld,x_h,y_1,\ld,y_k\rangle+\nu$ is in
$\fx$ for $j=1,\ld,m$. So we get the desired cones in $\fy$ and
$r_Y(P)$ is contractible.
\end{dimo}
\noindent\emph{Remark:} from the proofs of lemma \ref{a} and \ref{c}
we see  that:

\noindent for $P$ of type $a$,\ \  $r_X(P')=r_Y(P)+c\delta$, where 
$c=\min_{\,i} |r_Y(P)\cdot V(v_i)|$;\\
\noindent for $P$ of type $b$,\ \  $r_X(P')=r_Y(P)-\delta$; \\
\noindent for $P$ of type $c$,\ \  $r_X(P)=r_Y(P)$ and if 
$r_X(P')$ is contractible, 
 $r_X(P')=r_Y(P)-\delta$.

\end{document}